# GALOIS REPRESENTATIONS WITH CONJECTURAL CONNECTIONS TO ARITHMETIC COHOMOLOGY

AVNER ASH, DARRIN DOUD, AND DAVID POLLACK

Abstract. In this paper we extend a conjecture of Ash and Sinnott relating niveau one Galois representation to the mod $p$ cohomology of congruence subgroups of $\mathrm{SL}_n(\mathbb{Z})$ to include Galois representations of higher niveau. We then present computational evidence for our conjecture in the case $n = 3$ in the form of three-dimensional Galois representations which appear to correspond to cohomology eigenclasses as predicted by the conjecture. Our examples include Galois representations with nontrivial weight and level, as well as irreducible three-dimensional representations which are in no obvious way related to lower dimensional representations. In addition, we prove that certain symmetric square representations are actually attached to cohomology eigenclasses predicted by the conjecture.

## 1. Introduction

In [21], Serre published his conjecture (which had existed in some form since 1973) relating continuous odd absolutely irreducible Galois representations $\rho : G_{\mathbb{Q}} \to \mathrm{GL}_2(\bar{\mathbb{F}}_p)$ with the mod $p$ reductions of modular forms. He not only conjectured that a relationship existed, but also gave precise formulae describing where to find the predicted modular forms.

In [4], Ash and Sinnott present a conjecture giving a relationship between odd niveau 1 Galois representations of arbitrary dimension $n$, and certain cohomology groups of congruence subgroups of $\mathrm{GL}_n(\mathbb{Z})$. In the two-dimensional case, this conjecture is closely related to Serre's conjecture. Ash and Sinnott present computational evidence for their conjecture in certain three-dimensional cases, primarily in the case of three-dimensional level 1 reducible representations. In this paper we present additional computational evidence for the conjecture, including cases with nontrivial weight, level and nebentype. We also expand the conjecture to include representations of higher niveau and present computational evidence for this generalization. The representations in section 7.2 are particularly interesting—they are the first examples in which a cohomology eigenclass seems to correspond to a native three-dimensional Galois representation (i.e. a Galois representation which is in no obvious way related to a one-dimensional or two-dimensional Galois representation).

There is no problem in finding many Galois representations to which the conjecture that we make applies. The challenge is in finding Galois representations for

*Date*: January 26, 2001.

1991 *Mathematics Subject Classification.* 11F75,11F80.

The first author thanks the NSF and the NSA for their support. This manuscript is submitted for publication with the understanding that the United States government is authorized to reproduce and distribute reprints. The second author was supported by an NSF Postdoctoral Research Fellowship while doing this research.

which the predicted weights and levels allow computation of the associated cohomology classes and their Hecke eigenvalues. Our verifications of the conjecture involve finding representations which have fairly small weight and level, and computing the predicted cohomology groups and the action of Hecke operators on these groups for primes up to 47. We then compare the Hecke eigenvalues with the coefficients of the characteristic polynomials of the images of Frobenius, and if they match we claim to have evidence for the conjecture. We present many examples of Galois representations with weight and level small enough for us to work with, resulting in over 200 predictions (counting each weight associated to a Galois representation by Conjecture 3.1 separately). These examples are summarized in Tables 1 through 8, in which we describe Galois representations and give predicted weights, levels, and characters. For all the examples listed in the tables we have computed the homology groups (which are naturally dual to the cohomology groups), and in all cases an eigenclass with the correct eigenvalues up to $\ell = 47$ ($\ell = 3$ in Table 1) did exist in the predicted weight, level and character. Our examples include cases with niveau 1, 2, and 3, as well as wildly ramified niveau 1 representations. We also call attention to Theorem 4.1 and the examples which follow it, in which the theory of symmetric squares is used to *prove* a prediction of Conjecture 3.1 for certain irreducible three-dimensional Galois representations.

## 2. Definitions

Let $p$ be a prime number and $\bar{\mathbb{F}}_p$ an algebraic closure of the finite field $\mathbb{F}_p$ with $p$ elements. By a Galois representation, we mean a continuous representation of the absolute Galois group $G_{\mathbb{Q}}$ of $\mathbb{Q}$ to a matrix group $\mathrm{GL}_n(\bar{\mathbb{F}}_p)$. The representations with which we work in this paper will always be, in addition, semisimple. We say that a Galois representation is odd if the image of complex conjugation is a nonscalar matrix, and that it is even if the image of complex conjugation is a scalar.

For a given prime $q$ we will denote a decomposition group at $q$ in $G_{\mathbb{Q}}$ by $G_q$. This decomposition group then has a filtration by ramification subgroups $G_{q,i}$, with the whole inertia group above $q$ equal to $G_{q,0}$. We will often denote the inertia group $G_{p,0}$ at $p$ by $I_p$. We fix a Frobenius element $\mathrm{Frob}_q$ for each $q$, and we fix a complex conjugation $\mathrm{Frob}_\infty$.

We will denote the fundamental characters of niveau $n$ in characteristic $p$ [18] by $\psi_{n,d}$, $d = 1, \dots, n$, and note that they are all Galois conjugates (over $\mathbb{F}_p$) of $\psi_{n,1}$. In many cases we will be interested in working with fundamental characters of niveau 2 and 3, so for brevity we will let $\psi = \psi_{2,1}$, $\psi' = \psi_{2,2}$, and $\theta = \psi_{3,1}$, $\theta' = \psi_{3,2}$, and $\theta'' = \psi_{3,3}$. Note that the cyclotomic character $\omega$ is equal to $\psi_{1,1}$.

### 2.1. Hecke Operators.

Let $\Gamma_0(N)$ be the subgroup of matrices in $\mathrm{SL}_n(\mathbb{Z})$ whose first row is congruent to $(*, 0, \dots, 0)$ modulo $N$. Define $S_N$ to be the subsemigroup of integral matrices in $\mathrm{GL}_n(\mathbb{Q})$ satisfying the same congruence condition, and having positive determinant relatively prime to $N$.

Let $\mathcal{H}(N)$ denote the $\bar{\mathbb{F}}_p$-algebra of double cosets $\Gamma_0(N)S_N\Gamma_0(N)$. Then $\mathcal{H}(N)$ is a commutative algebra which acts on cohomology and homology of $\Gamma_0(N)$ with coefficients in any $\bar{\mathbb{F}}_p[S_N]$-module. When a double coset is acting on cohomology or homology we call it a Hecke operator. Clearly, $\mathcal{H}(N)$ contains all double cosets of the form $\Gamma_0(N)D(\ell,k)\Gamma_0(N)$, where $\ell$ is a prime not dividing $N$, $0 \le k \le n$, and

$$D(\ell,k) = \begin{pmatrix} 1 & & & & & \\ & \ddots & & & & \\ & & 1 & & & \\ & & & \ell & & \\ & & & & \ddots & \\ & & & & & \ell \end{pmatrix}$$

is the diagonal matrix with the first $n-k$ diagonal entries equal to 1 and the last $k$ diagonal entries equal to $\ell$. When we consider the double coset generated by $D(\ell,k)$ as a Hecke operator we call it $T(\ell,k)$.

**Definition 2.1.** Let $V$ be an $\mathcal{H}(pN)$-module, and suppose that $v \in V$ is a simultaneous eigenvector for all $T(\ell,k)$ and that $T(\ell,k)v = a(\ell,k)v$ with $a(\ell,k) \in \bar{\mathbb{F}}_p$ for all $\ell \not| \, pN$ prime and all $0 \le k \le n$. If

$$\rho : G_{\mathbb{Q}} \to \mathrm{GL}_n(\bar{\mathbb{F}}_p)$$

is a representation unramified outside $pN$ and

$$\sum_{k=0}^{n} (-1)^k \ell^{k(k-1)/2} a(\ell,k) X^k = \det(I - \rho(\mathrm{Frob}_\ell)X)$$

for all $\ell \not| \, pN$, then we say that $\rho$ is attached to $v$ or that $v$ corresponds to $\rho$.

2.2. **Level and Nebentype.** Let

$$\rho : G_{\mathbb{Q}} \to \mathrm{GL}_n(\bar{\mathbb{F}}_p)$$

be a continuous representation. We will define a level and nebentype associated to $\rho$ exactly as Serre does in [21].

For a fixed prime $q \ne p$, and $i \ge 0$, let $g_i = |\rho(G_{q,i})|$. Note that this is a finite integer, since by continuity the image of $\rho$ must be finite. Let $M = \bar{\mathbb{F}}_p^n$ be acted on by $G_{\mathbb{Q}}$ via $\rho$ in the natural way, and define

$$n_q = \sum_{i=0}^{\infty} \frac{g_i}{g_0} \dim M/M^{G_{q,i}}.$$

The sum defining $n_q$ is then a finite sum, since eventually the $G_{q,i}$ are trivial.

**Definition 2.2.** With $\rho$ as above, define the level

$$N(\rho) = \prod_{q \ne p} q^{n_q}.$$

Note that this product is also finite, since $\rho$ is ramified at only finitely many primes, and $n_q$ is 0 if $\rho$ is unramified at $q$.

In order to define the nebentype character, we again proceed exactly as Serre does in [21]. We factor $\det \rho = \epsilon\omega^k$, where $\omega$ is the cyclotomic character modulo $p$, and $\epsilon$ is a character $G_{\mathbb{Q}} \to \bar{\mathbb{F}}_p$ unramified at $p$. By class field theory, we may then consider $\epsilon$ as a Dirichlet character

$$\epsilon : (\mathbb{Z}/N(\rho)\mathbb{Z})^* \to \bar{\mathbb{F}}_p.$$

We then pull back the definition of $\epsilon$ to $S_N$ by defining $\epsilon$ to be the composite character

$$S_N \to (\mathbb{Z}/N(\rho)\mathbb{Z})^* \to \bar{\mathbb{F}}_p,$$

where the first map takes a matrix in $S_N$ to its $(1,1)$ entry, and define $\mathbb{F}_\epsilon$ to be the one-dimensional space $\bar{\mathbb{F}}_p$ with the action of $S_N$ given by $\epsilon$.

We now define

$$V(\epsilon) = V \otimes \mathbb{F}_\epsilon.$$

Letting $\Gamma_0(N)$ act on $V$ by reduction modulo $p$, we see that $V(\epsilon)$ is a $\Gamma_0(N)$-module. In addition, letting $S_{pN}$ act via $\epsilon$ on $\bar{\mathbb{F}}_p$, we see that $V(\epsilon)$ is also an $S_{pN}$-module.

In specifying the nebentype, we will often refer to the unique quadratic character modulo $p$ ramified only at a prime $q > 3$, and will denote this character by

$$\epsilon_q : G_{\mathbb{Q}} \to \mathbb{F}_p.$$

We will also refer to the character $\epsilon_4$, which is ramified only at 2 and cuts out the field $\mathbb{Q}(\sqrt{-1})$.

2.3. **Irreducible $\mathbf{GL}_n(\mathbb{F}_p)$-modules.** The natural generalization of the weight in Serre's conjecture is an irreducible $\mathrm{GL}_n(\mathbb{F}_p)$-module. To see this, we note that the Eichler-Shimura theorem [22] relates the space of modular forms of weight $k$ to cohomology with coefficients in

$$\mathrm{Sym}^g(\mathbb{C}^2)$$

with $g = k-2$. Hence, an eigenform $f$ of level $N$, nebentype $\epsilon$, and weight $k$ gives rise to a collection of Hecke eigenvalues which, when reduced modulo $p$, also occur in

$$H^1(\Gamma_0(N), V_g(\epsilon)),$$

where $V_g \cong \mathrm{Sym}^g(\bar{\mathbb{F}}_p^2)$ is the space of two-variable homogeneous polynomials of degree g over $\bar{\mathbb{F}}_p$ with the natural action of $\mathrm{SL}_2(\bar{\mathbb{F}}_p)$. Ash and Stevens have shown [5] that any system of Hecke eigenvalues occurring in cohomology with coefficients in some $\mathrm{GL}_n(\mathbb{F}_p)$-module also occurs in the cohomology with coefficients in at least one irreducible $\mathrm{GL}_n(\mathbb{F}_p)$-module occurring in a composition series of the original module. Hence, there is some irreducible $\mathrm{GL}_n(\mathbb{F}_p)$-module $W$, such that the system of eigenvalues coming from $f$ also occurs in $H^1(\Gamma_0(N), W(\epsilon))$. Given this fact, it is natural to ask which irreducible component(s) gives rise to the system of eigenvalues.

We may parameterize the complete set of irreducible $\mathrm{GL}_n(\mathbb{F}_p)$-modules as in [10].

**Definition 2.3.** We say that an $n$-tuple of integers $(b_1, \dots, b_n)$ is $p$-restricted if

$$0 \le b_i - b_{i+1} \le p-1, \quad 1 \le i \le n-1,$$

and

$$b_n < p-1.$$

**Proposition 2.4.** *The set of irreducible $GL_n(\mathbb{F}_p)$-modules is in one-to-one correspondence with the collection of all $p$-restricted $n$-tuples.*

The one-to-one correspondence in this proposition may be described explicitly as follows: the module $F(b_1, \dots, b_n)$ corresponding to the $p$-restricted $n$-tuple $(b_1, \dots, b_n)$ is the unique simple submodule of the dual Weyl module $W(b_1, \dots, b_n)$ with highest weight $(b_1, \dots, b_n)$. Theorem 8.1 gives an explicit model for the module $F(b_1, b_2, b_3)$ in the case $n = 3$, but for larger $n$ no general computational models are known to the authors.

In dealing with Galois representations it will often become necessary to associate an irreducible module to an $n$-tuple which is not $p$-restricted. We do this via the following definition:

**Definition 2.5.** Let $(a_1, \dots, a_n)$ be any $n$-tuple of integers. Define

$$F(a_1, \dots, a_n)' = F(b_1, \dots, b_n),$$

where $(b_1, \dots, b_n)$ is a $p$-restricted $n$-tuple for which

$$a_i \equiv b_i \pmod{p-1}.$$

We note that in certain cases (namely when some $a_i \equiv a_{i+1} \pmod{p-1}$) the module $F(a_1, \dots, a_n)'$ may not be well defined. In this case, we interpret any statement concerning $F(a_1, \dots, a_n)'$ to mean that the statement is true for some choice of $F(b_1, \dots, b_n)$ as in the definition. For example, if $p = 5$, a statement concerning $F(1,0,0)'$ is true if the statement is true for either $F(1,0,0)$ or $F(5,4,0)$ (or both). When dealing with modules defined by the prime notation, we will say that a module $F(a_1, \dots, a_n)'$ is determined unambiguously if there is a unique $p$-restricted sequence congruent to $(a_1, \dots, a_n)$ modulo $p-1$.

2.4. **The strict parity condition.** We modify slightly the statement of the strict parity condition in [4] for ease of exposition, but our formulation is logically equivalent to that in [4].

**Definition 2.6.** Let $V = \bar{\mathbb{F}}_p^n$ be an $n$-dimensional space, with the standard action of $\mathrm{GL}_n(\bar{\mathbb{F}}_p)$. A Levi subgroup $L$ of $\mathrm{GL}_n(\bar{\mathbb{F}}_p)$ is the simultaneous stabilizer of a collection $D_1, \dots, D_k$ of subspaces such that $V = \bigoplus_i D_i$. If each $D_i$ has a basis consisting of standard basis vectors for $V$, then $L$ is called a standard Levi subgroup.

*Example* 2.7. The standard Levi subgroups of $\mathrm{GL}_2(\bar{\mathbb{F}}_p)$ are the subgroup of diagonal matrices and the whole of $\mathrm{GL}_2(\bar{\mathbb{F}}_p)$.

*Example* 2.8. The standard Levi subgroups of $\mathrm{GL}_3(\bar{\mathbb{F}}_p)$ are the subgroup of diagonal matrices, the whole of $\mathrm{GL}_3(\bar{\mathbb{F}}_p)$, and the three subgroups

$$\begin{pmatrix} * & 0 & 0 \\ 0 & * & * \\ 0 & * & * \end{pmatrix}, \quad \begin{pmatrix} * & 0 & * \\ 0 & * & 0 \\ * & 0 & * \end{pmatrix}, \quad \begin{pmatrix} * & * & 0 \\ * & * & 0 \\ 0 & 0 & * \end{pmatrix}.$$

**Definition 2.9.** Let $\rho : G_{\mathbb{Q}} \to \mathrm{GL}_n(\bar{\mathbb{F}}_p)$ be a continuous representation. A standard Levi subgroup $L$ of $\mathrm{GL}_n(\bar{\mathbb{F}}_p)$ is said to be $\rho$-minimal if $L$ is minimal among all standard Levi subgroups which contain some conjugate of the image of $\rho$.

**Definition 2.10.** A semisimple continuous representation $\rho : G_{\mathbb{Q}} \to \mathrm{GL}_n(\bar{\mathbb{F}}_p)$ satisfies the strict parity condition with Levi subgroup $L$ if it has the following properties:

1. The image of $\rho$ lies inside a $\rho$-minimal standard Levi subgroup $L$,
2. The image of complex conjugation is conjugate inside $L$ to a matrix

$$\pm \begin{pmatrix} 1 & & \\ & -1 & \\ & & \ddots \end{pmatrix}$$

with strictly alternating signs on the diagonal.

*Example* 2.11. Any odd irreducible two-dimensional (resp. three-dimensional) representation satisfies strict parity, with $L = \mathrm{GL}_2(\bar{\mathbb{F}}_p)$ (resp. $L = \mathrm{GL}_3(\bar{\mathbb{F}}_p)$).

*Example* 2.12. Let $\rho$ be the direct sum of a two-dimensional odd irreducible representation and a one-dimensional representation, with image contained inside

$$L = \begin{pmatrix} * & 0 & 0 \\ 0 & * & * \\ 0 & * & * \end{pmatrix} \quad \textit{or} \quad L = \begin{pmatrix} * & * & 0 \\ * & * & 0 \\ 0 & 0 & * \end{pmatrix}.$$

Then $\rho$ satisfies the strict parity condition, with Levi subgroup $L$.

*Example* 2.13. Let $\rho$ be the direct sum of a two-dimensional even irreducible representation and a one-dimensional representation, with the image of $\rho$ contained inside

$$L = \begin{pmatrix} * & & * \\ & * & \\ * & & * \end{pmatrix}.$$

Then $\rho$ satisfies strict parity with this Levi subgroup exactly when $\rho$ is odd.

*Remark* 2.14. Note that any odd three-dimensional Galois representation is conjugate to a representation which satisfies the strict parity condition for some standard Levi subgroup $L$. More generally, if $/rho$ is an $n$-dimensional representation where the number of $+1$ eigenvalues and the number of $-1$ eigenvalues of complex conjugation differ by at most one, then $\rho$ satisfies the strict parity condition for some standard Levi subgroup L.

**Definition 2.15.** If $\rho : G \to \mathrm{GL}_n(\bar{\mathbb{F}}_p)$ lands inside a Levi subgroup $L$, and $\sigma : G \to \mathrm{GL}_n(\bar{\mathbb{F}}_p)$ is another representation of $G$, we say that

$$\rho \sim_L \sigma$$

if there is a matrix $M \in L$ such that

$$M\rho(g)M^{-1} = \sigma(g)$$

for all $g \in G$. If $L = \mathrm{GL}_n(\bar{\mathbb{F}}_p)$, then we may write

$$\rho \sim \sigma.$$

2.5. **Weights.** We now begin to predict the weights (or irreducible modules) for which we expect to find cohomology eigenclasses with $\rho$ attached. Following the example of Serre's conjecture, we expect these weights to be determined by the restriction of $\rho$ to a decomposition group at $p$, so we are interested in studying representations of the decomposition group $G_p$. For convenience, we will denote the inertia group $G_{p,0}$ by $I_p$ and the wild ramification group $G_{p,1}$ by $I_w$. We begin by considering simple representations of $G_p$.

**Lemma 2.16.** *Let $V$ be a simple $n$-dimensional $\bar{\mathbb{F}}_p[G_p]$-module, with the action of $G_p$ given by a representation $\rho : G_p \to GL(V)$. Then we may choose a basis for $V$ such that*

$$\rho|_{I_p} = \begin{pmatrix} \varphi_1 & & \\ & \ddots & \\ & & \varphi_n \end{pmatrix},$$

*with the characters $\varphi_1, \dots, \varphi_n$ equal to some permutation of $\psi_{n,1}^m, \dots, \psi_{n,n}^m$ for some $m \in \mathbb{Z}$.*

*Proof.* This proof is almost identical to the proof in [13] for two-dimensional representations. We first note that $\rho$ has finite image, so that we may actually realize it over a finite extension of $\mathbb{F}_p$. Hence, we may find an $\mathbb{F}_{p^m}[G_p]$-module $V'$, such that $V = V' \otimes \bar{\mathbb{F}}_p$. We note that $I_w$ must act trivially on $V'$, since the invariants ${V'}^{I_w}$ are a $G_p$-submodule of the simple module $V'$, and also nontrivial (since the image of $I_w$ under $\rho$ is a $p$-group). Hence, we may diagonalize $\rho|_{I_p}$. Since the Frobenius acts on the tame inertia as $p$th powers, we see that the set of diagonal characters must be stable under taking $p$th powers. Finally, since $V$ is simple, the Frobenius must permute the diagonal characters transitively, resulting in the characterization given above. □

*Remark* 2.17. Note that for a given $V$, Lemma 2.16 yields $n$ distinct values of $m$ modulo $(p^n-1)$. If $m_0$ is one of them, the others are congruent to $pm_0, p^2m_0, \dots, p^{n-1}m_0$ modulo $(p^n - 1)$.

**Definition 2.18.** Let $V$ be a simple $G_p$-module, diagonalized as in Lemma 2.16 with some choice of exponent $m$. Write $m$ as

$$m = a_1 + a_2 p + \cdots + a_n p^{n-1},$$

with $0 \le a_i - a_n \le p-1$ for all $i$. Suppose that $(b_1, \dots, b_n)$ satisfies $b_i \ge b_{i+1}$ for all $i < n$ and is obtained by permuting the entries of $(a_1, \dots, a_n)$. Then we say that $(b_1, \dots, b_n)$ is derived from $\rho$. If the action of $G_p$ on $V$ is given by a representation $\rho$ we say that the $n$-tuple is derived from $\rho$.

*Remark* 2.19. Note that not all values of $m$ have an expansion of the form given here. For example, if $p = 5$, $n = 3$, $m = 30$, there is no expansion satisfying the above properties. It is a simple exercise to see that every simple module has at least one derived $n$-tuple, and that a given value of $m$ yields a unique $n$-tuple, if it yields any. Hence a simple $n$-dimensional $G_p$-module may have at most $n$ $n$-tuples derived from it, but can have fewer.

Now let $V$ be any $n$-dimensional $G_p$-module, with the action of $G_p$ given by $\rho : G_p \to \mathrm{GL}(V)$. We may find a composition series

$$\{0\} = V_0 \subset V_1 \subset \cdots \subset V_k = V.$$

Let each composition factor $V_i/V_{i-1}$ have dimension $d_i$, and set $d_0 = 0$.

By diagonalizing $\rho$ on each simple composition factor, we may find a basis $(e_1, \dots, e_n)$ of $V$ such that $\rho$ is upper triangular, with diagonal characters

$$(\varphi_{1,1}, \dots, \varphi_{1,d_1}, \varphi_{2,1}, \dots, \varphi_{2,d_2}, \dots, \varphi_{k,1}, \dots, \varphi_{k,d_k}),$$

where the first $d_1$ characters come from the action on $V_1/V_0$, the next $d_2$ from the action on $V_2/V_1$, etc. For each composition factor, choose $m_i$ so that for some $j$, $\psi_{d_i,1}^{m_i} = \varphi_{i,j}$, and so that $m_i$ yields a $d_i$-tuple derived from $V_i/V_{i-1}$. Concatenating these $d_i$-tuples gives us an $n$-tuple $(a_1, \dots, a_n)$.

We wish to preserve the order of the integers in our $n$-tuple which come from an individual composition factor, so we make the following definition.

**Definition 2.20.** A permutation $\sigma$ of the integers $\{1, \dots, n\}$ is compatible with the filtration

$$0 = V_0 \subset V_1 \subset \ldots \subset V_k = V$$

given above, if for $0 \le s < k$, and $a, b \in [1 + \sum_{j=0}^{s} d_j, d_{s+1} + \sum_{j=0}^{s} d_j]$ with $a < b$, we have $\sigma(a) < \sigma(b)$.

**Definition 2.21.** Let $V$ be an $n$-dimensional $G_p$-module with chosen basis $\{e_1, \dots, e_n\}$ with respect to which the action of $G_p$ is upper triangularized, and let $(a_1, \dots, a_n)$ be an $n$-tuple obtained as above. If $\sigma$ is a permutation of the integers $\{1, \dots, n\}$ compatible with the filtration above and such that the action of $G_p$ with respect to the ordered basis $\{e_{\sigma(1)}, \dots, e_{\sigma(n)}\}$ remains upper triangular, then we will say that the $n$-tuple $(a_{\sigma(1)}, \dots, a_{\sigma(n)})$ is derived from $V$.

*Remark* 2.22. Note that there is at least one (and possibly more) $n$-tuple derived from $V$, namely the original $n$-tuple $(a_1, \dots, a_n)$. In addition, even the choice of this original $n$-tuple is not unique, so that there will usually be many $n$-tuples derived from a given $V$.

**Definition 2.23.** Let $\rho : G_{\mathbb{Q}} \to GL_n(\bar{\mathbb{F}}_p)$ be a semisimple continuous representation, conjugated to land in a $\rho$-minimal standard Levi subgroup $L$. Let $D_1, \dots, D_k$ be the subspaces of $\bar{\mathbb{F}}_p^n$ given in the definition of $L$. Then we have representations $\rho_i : G_{\mathbb{Q}} \to GL(D_i)$, which make each $D_i$ into a $G_{\mathbb{Q}}$-module. Let $G_p$ be a decomposition group above $p$, and consider each $D_i$ as a $G_p$-module. Let $d_i = \dim D_i$, and let $(a_1, \dots, a_{d_i})$ be a $d_i$-tuple derived from $D_i$ as above. If the standard basis elements of $\bar{\mathbb{F}}_p^n$ which span $D_i$ are $e_{j_r}$ $(1 \le r \le d_i)$ with $j_r < j_s$ for $r < s$, then set $b_{j_r} = a_r$ for $r = 1, \dots, d_i$. Doing this for each $D_i$ produces an $n$-tuple $(b_1, \dots, b_n)$. Such an $n$-tuple will be said to be derived from $\rho$, with Levi subgroup $L$.

*Remark* 2.24. Note that the above discussion may (in many cases) be summarized more informally as follows: given a representation $\rho : G_{\mathbb{Q}} \to \mathrm{GL}_n(\bar{\mathbb{F}}_p)$ which lands inside a $\rho$-minimal standard Levi subgroup $L$, we may upper triangularize its restriction to inertia by conjugating by an element of $L$. This will give a sequence of characters of the tame inertia group on the diagonal. Group these characters together into "niveau $d$ collections" (a niveau $d$ collection is set of $d$ characters, each a power of a different fundamental character of niveau $d$ with the same exponent $m$, and all appearing in the same "Levi Block"). For a given niveau $d$ collection, write the exponent $m$ as $a_1 + a_2 p + \cdots + a_d p^{d-1}$, with $0 \le a_i - a_n \le p - 1$ for all $i$, and let $(b_1, \dots, b_d)$ be the ordered (decreasing) $d$-tuple with the same components as $(a_1, \dots, a_d)$. Then construct an $n$-tuple $(c_1, \dots, c_n)$ as follows: if the $i$th character in the niveau $d$ collection is in the $k$th diagonal position in the image of $\rho$, set $c_k = b_i$ (note that the order of the $b_i$ should be preserved in the $n$-tuple). This procedure gives the same derived $n$-tuples as above, except when there is a combination of wild ramification and multiple niveau $d$ collections containing the same characters, in which case the more complicated procedure described above is needed.

## 3. Conjecture

**Conjecture 3.1.** *Let $\rho : G_{\mathbb{Q}} \to GL_n(\bar{\mathbb{F}}_p)$ be a continuous semisimple Galois representation. Suppose that $\rho$ satisfies the strict parity condition with Levi subgroup $L$. Let $(a_1, \dots, a_n)$ be an $n$-tuple derived from $\rho$ with the Levi subgroup $L$, and define $V = F(a_1 - (n-1), a_2 - (n-2), \dots, a_n - 0)'$. Further, let $N = N(\rho)$ be the level of $\rho$ and let $\epsilon = \epsilon(\rho)$ be the nebentype character of $\rho$. Then $\rho$ is attached to a cohomology eigenclass in*

$$H^*(\Gamma_0(N), V(\epsilon)).$$

*Remark* 3.2. We note that in the case of two-dimensional Galois representations, we may take the $*$ to be 0 or 1, and in fact for irreducible two-dimensional representations we may take $*$ to be 1.

In the case of three-dimensional Galois representations we may take $*$ to be at most 3, and for irreducible Galois representations (or sums of an even two-dimensional representation with a one-dimensional representation) we may take $*$ to be equal to 3, as explained in [4]. As mentioned previously, any odd two-dimensional or three-dimensional representation is conjugate to a representation which satisfies strict parity for some standard Levi subgroup $L$.

In our computations, we test the conjecture for three-dimensional representations by computing $H^3$. In cases where $\rho$ is the sum of three characters or the sum of an odd two-dimensional representation and a character we are thus actually testing a stronger assertion than Conjecture 3.1, namely that the cohomology class exists in $H^3$. See for example Tables 3 and 7. We did not test any $\rho$ which are sums of three characters in this paper, but several examples of such may be found in [3] and [1]. In addition, we do not present computational examples for $p = 2$, as this would involve rewriting portions of our computer programs. In addition, for $p = 2$ and $p = 3$, our computational techniques (based on those in [1]) do not always compute the whole of $H^3$. Nevertheless, we have no reason to doubt our conjecture for these primes. In particular, problems with the weight and nebentype that occur when $p = 2$ or $p = 3$ for Serre's original conjecture involving classical modular forms modulo $p$ should not occur for our conjecture which involves mod $p$ cohomology.

*Remark* 3.3. Note that Conjecture 3.1 applies to Galois representations of arbitrary dimension, but we have no computational evidence for dimension higher than 3. Forthcoming work of the first author with P. Gunnells and M. McConnell touches on the case of certain four-dimensional representations.

*Remark* 3.4. Note that the conjecture makes no claim of predicting all possible weights which yield an eigenclass with $\rho$ attached. In fact we have three types of computational examples in which additional weights (not predicted by the conjecture) do yield eigenclasses which appear to have $\rho$ attached.

The first type of additional weight occurs if $\rho$ is attached to a "quasi-cuspidal" eigenclass (for instance if $\rho$ is either irreducible or reducible as a sum of an even two-dimensional representation and a character). In this case, for certain weights, we may define an extra weight as follows:

**Definition 3.5.** Let $F(a,b,c)$ be an irreducible $\mathrm{GL}_n(\mathbb{F}_p)$-module, with $a-c < p-2$. Then we may define

$$M = F(d,e,f) = \begin{cases} F(p-2+c, b, a-(p-2)) & \text{if } a \geq p-2 \\ F(2(p-2)+c+1, b+(p-1), a+1) & \text{if } a < p-2 \end{cases}.$$

Then we say that $M$ is the extra weight associated to $F(a,b,c)$.

Applying [10, Proposition 2.11] it is easy to see that if $F(d,e,f)$ is the extra weight associated to $F(a,b,c)$, there is an exact sequence

$$0 \to F(d,e,f) \to W(d,e,f) \to F(a,b,c) \to 0.$$

Now, suppose that $\rho$ is attached to a quasi-cuspidal cohomology eigenclass in weight $F(a,b,c)$. Examining the long exact homology sequence associated to this short exact sequence, we find that a quasi-cuspidal eigenclass $\alpha$ in $H_3(\Gamma_0(N), F(a,b,c)(\epsilon))$ (in particular, any eigenclass corresponding to an irreducible Galois representation) is either the image of an eigenclass in $H_3(\Gamma_0(N), W(d,e,f)(\epsilon))$, or has nonzero image $\beta$ in $H_2(\Gamma_0(N), F(d,e,f)(\epsilon))$. In the second case, $\beta$ is an eigenclass, and using Theorem 3.10 and Lefschetz duality we find that there is an eigenclass $\gamma$ in $H_3(\Gamma_0(N), F(d,e,f)(\epsilon))$ which has the same eigenvalues as $\alpha$. Hence, for each quasi-cuspidal eigenclass in an appropriate weight, there are two possibilities: either the eigenclass lifts to the dual Weyl module, or the eigenclass gives rise to another eigenclass with the same eigenvalues in the extra weight. Our experimental evidence supports the hypothesis that in all such cases a quasi-cuspidal eigenclass gives rise to another eigenclass with the same eigenvalues in the extra weight.

The second class of additional weights which we have observed consists of certain weights which would be predicted by our conjecture if we eliminated the strict parity condition. These additional weights have only been observed for representations $\rho$ which are either the sum of three characters or the sum of an odd two-dimensional representation and a character. These additional weights seem to occur fairly rarely and sporadically, and may be related to the occurrence of eigenclasses in $H^2$ which have $\rho$ attached. A full investigation of them would require new computational techniques, beyond those developed in this paper.

The third class of additional weights consists of extra weights associated to weights which would be predicted by our conjecture but for the strict parity condition. As in the second case, these additional weights only occur rarely, and only for reducible $\rho$.

Before beginning to present computational evidence for Conjecture 3.1, we begin by proving several facts about the conjecture.

**Theorem 3.6.** *If Conjecture 3.1 is true for a representation $\rho$, then it is true for the representation $\rho \otimes \omega^s$, where $\omega$ is the cyclotomic character modulo $p$.*

*Proof.* First note that twisting by $\omega^s$ does not affect the predicted level or nebentype in any way. Denote the level of $\rho$ by $N$ and the nebentype of $\rho$ by $\epsilon$.

If $\rho$ has niveau 1, then this is just Proposition 2.6 of [4].

For higher niveau representations, we note that twisting by $\omega^s$ changes the value of $m$ coming from a niveau $d$ character by $s(1+p+\cdots+p^{d-1})$, hence it changes all the values of $a_i$ arising from $m$ by $s$. Following this change through the permutations involved in deriving an $n$-tuple, we find that twisting a representation $\rho$ by $\omega^s$ adds $s$ to each element of a derived $n$-tuple. This change is then reflected in the predicted weight, and we have that the set of predicted weights for $\rho \otimes \omega^s$ is precisely the set of twists by $\det^s$ of the predicted weights of $\rho$.

Finally, if an eigenclass $v$ shows up in weight $V$, and has $\rho$ attached, then we may consider $v$ as lying in cohomology with weight $V \otimes \det^s$, and see easily (as in [4]) that in this new cohomology group $v$ has $\rho \otimes \omega^s$ attached. Hence, if $\rho$ is attached to a cohomology class in each of the weights predicted by Conjecture 3.1, then $\rho \otimes \omega^s$ will be as well. ∎

We now note that there is a correspondence between systems of Hecke eigenvalues arising from modular forms and systems of eigenvalues arising from arithmetic cohomology in characteristic $p$, similar to that given by the Eichler-Shimura isomorphism in characteristic 0. In particular, we note that by [6, Proposition 2.5], for $p > 3$, any system of Hecke eigenvalues comes from the mod $p$ reduction of an eigenform of level $N$, nebentype $\epsilon$ and weight $k = g + 2$ if and only if is comes from a Hecke eigenclass in $H^1(\Gamma_0(N), V_g(\mathbb{F}_p)(\epsilon))$, where $V_g(\mathbb{F}_p)$ is the $g$th symmetric power of the standard representation of $\mathrm{GL}_2(\mathbb{F}_p)$.

**Theorem 3.7.** *If $p > 3$, Serre's conjecture implies Conjecture 3.1 for $n = 2$.*

*Proof.* For a complete description of Serre's conjecture, including Serre's prediction of the weight, see [21] or [13].

There are two cases, where $\rho$ is niveau 1 or niveau 2. In either case, we note that the level and nebentype predicted by Serre's conjecture are identical to those predicted by Conjecture 3.1, so that we need only deal with the weight.

Suppose that $\rho : G_{\mathbb{Q}} \to GL_2(\bar{\mathbb{F}}_p)$ is odd, semisimple, and has niveau 1. If $\rho$ is reducible, Conjecture 3.1 is true [4, Proposition 2.7], so we may assume that $\rho$ is irreducible. If $\rho$ is tamely ramified, we have that

$$\rho|_{I_p} \sim \begin{pmatrix} \omega^{a_1} & \\ & \omega^{a_2} \end{pmatrix},$$

with $0 \le a_1, a_2 \le p-2$. Conjecture 3.1 predicts a weight of $F(a_1 - 1, a_2)'$.

If $a_2 < a_1$, then

$$\rho \otimes \omega^{-a_2}|_{I_p} \sim \begin{pmatrix} \omega^{a_1-a_2} & \\ & \omega^0 \end{pmatrix},$$

and Serre's conjecture claims that $\rho \otimes \omega^{-a_2}$ corresponds to a modular form of weight $1 + a_1 - a_2$, or (via [6, Proposition 2.5]) that $\rho \otimes \omega^{-a_2}$ corresponds to a cohomology class with coefficients in $F(a_1 - a_2 - 1, 0)$. Twisting by $\omega^{a_2}$ (which corresponds to twisting the weight by $\det^{a_2}$), we find that $\rho$ corresponds to a cohomology class with coefficients in $F(a_1 - 1, a_2)$, exactly as predicted by Conjecture 3.1.

If $a_2 \ge a_1$, then

$$\rho \otimes \omega^{p-1-a_2}|_{I_p} \sim \begin{pmatrix} \omega^{p-1+a_1-a_2} & \\ & \omega^0 \end{pmatrix},$$

and by Serre's conjecture (together with [6, Proposition 2.5]), $\rho \otimes \omega^{p-1-a_2}$ corresponds to a cohomology class with coefficients in $F(p - 2 + a_1 - a_2, 0)$. Twisting by $\omega^{a_2}$, as before, we find that $\rho$ has weight $F(a_1 - 1 + (p - 1), a_2)$, exactly as predicted.

Now if $\rho$ is wildly ramified at $p$, then

$$\rho|_{I_p} \sim \begin{pmatrix} \omega^{\beta} & * \\ & \omega^{\alpha} \end{pmatrix},$$

with $0 \le \alpha \le p - 2$ and $1 \le \beta \le p - 1$, and Conjecture 3.1 predicts a weight of $F(\beta - 1, \alpha)'$. Before applying Serre's conjecture, we will twist $\rho$ by $\omega^{-\alpha}$, to obtain

$$\rho \otimes \omega^{-\alpha} \sim \begin{pmatrix} \omega^{\beta-\alpha} & * \\ & \omega^0 \end{pmatrix},$$

Applying Serre's conjecture to this representation, we find that it has weight

1. $1 + (\beta - \alpha)$ (or $F(\beta - \alpha - 1, 0)$) if $\beta > \alpha + 1$.
2. 2 (or $F(0, 0)$) if $\beta = \alpha + 1$ and $\rho \otimes \omega^{-\alpha}$ is *peu ramifiée.*
3. $p + 1$ (or $F(p - 1, 0)$) if $\beta = \alpha + 1$ and $\rho \otimes \omega^{-\alpha}$ is *très ramifiée.*
4. $1 + (\beta - \alpha) + (p - 1)$ (or $F(\beta - \alpha + (p - 1) - 1, 0)$ if $\beta \le \alpha$

Twisting each of these weights by $\det^{\alpha}$, we find that $\rho$ corresponds to a cohomology class in weight $F(\beta - 1, \alpha)'$ in every case (note that when $\beta - 1 \le \alpha$ we may add $p - 1$ to $\beta - 1$ to obtain a $p$-restricted pair).

This proves the theorem in case $\rho$ has niveau 1.

Suppose that

$$\rho|I_p \sim \begin{pmatrix} \psi^m & \\ & \psi'^m \end{pmatrix},$$

where $m = a + bp$, and $0 < a - b \leq p - 1$ (note that if $a = b$, we are really in niveau 1). For simplicity we use the fact that $\psi$ and $\psi'$ have order $p^2 - 1$ to reduce to the case where $0 < m < p^2 - 1$, so that $b < p - 1$. The weight predicted by Conjecture 3.1 is then $F(a-1, b)'$.

Now,

$$\rho \otimes \omega^{-b} \sim_L \begin{pmatrix} \psi^{a-b} & \\ & \psi'^{a-b} \end{pmatrix},$$

so that by Serre's conjecture $\rho \otimes \omega^{-b}$ corresponds to a cohomology class with coefficients in $F(a-b-1, 0)$. Twisting by $\omega^b$, we see that $\rho$ then corresponds to a cohomology class with coefficients in $F(a-1, b)$, exactly as predicted.

Hence, Serre's conjecture implies Conjecture 3.1, for $n = 2$. □

We now prove a partial converse to Theorem 3.7, which shows that in certain cases Conjecture 3.1 is actually equivalent to Serre's conjecture.

**Theorem 3.8.** *Assume Conjecture 3.1. Let $p > 3$, and let $\rho : G_{\mathbb{Q}} \to GL_2(\bar{\mathbb{F}}_p)$ be a semisimple continuous odd Galois representation. If each weight predicted by Conjecture 3.1 is defined unambiguously, then Serre's conjecture is true for $\rho$.*

*Proof.* We may clearly assume that $\rho : G_{\mathbb{Q}} \to \mathrm{GL}_2(\bar{\mathbb{F}}_p)$ is irreducible, since Serre's conjecture says nothing about reducible representations.

First, note that $\rho$ can not be attached to any class in $H^0$, since according to [2, Theorem 4.1.4], any class in $H^0$ is a twist of a "punctual" class, and a punctual class corresponds to a reducible representation by [2, Lemma 4.1.2].

Conjecture 3.1 implies that $\rho$ is attached to an eigenclass in $H^1(\Gamma_0(N), V(\epsilon))$, where $N$, $V$, and $\epsilon$ are as predicted in the conjecture. We note that the level and nebentype predicted by Conjecture 3.1 are exactly the same as those predicted by Serre's conjecture.

If $\rho$ is tamely ramified and has niveau 1, then we have

$$\rho|_{I_p} \sim \begin{pmatrix} \omega^a & \\ & \omega^b \end{pmatrix},$$

and we may further conjugate $\rho$ so that $0 \leq b \leq a < p - 1$. The weights predicted by Conjecture 3.1 are then $F(a-1, b)'$, and (permuting the diagonal characters) $F(b-1, a)'$. These will be defined unambiguously exactly when $a \neq b + 1$. For $a > b + 1$, we have that $F(a-1, b)$ embeds in $F(a-1+bp, 0)$, so conjecture 3.1 predicts that $\rho$ is attached to a cohomology eigenclass in weight $F(a-1+bp, 0)$, since any system of eigenvalues occurring in a submodule will occur in the containing module [4]. This implies (by [6, Proposition 2.5]) that $\rho$ is attached to an eigenform of weight $1+a+bp$, which is exactly the weight predicted by Serre's conjecture. For $a = b = 0$, the predicted weights for Conjecture 3.1 and Serre's conjecture are both $F(p-2, 0)$. For $a = b \neq 0$, Conjecture 3.1 predicts a weight of $F(a-1+p-1, b)$, while Serre's conjecture predicts a weight of $F(b-1+pa, 0)$. Using [10, Table 1] (specifically the last line, as $b \neq 0$), we see that $F(a-1+p-1, b)$ is a subquotient of $F(b-1+pa, 0)$. Hence we are finished if we can show that the system of eigenvalues corresponding to $\rho$ in weight $F(a-1+p-1, b)$ also shows up in weight $F(b-1+pa, 0)$. Lemma 3.9 shows that systems of eigenvalues of eigenclasses which are not twists of punctual classes are inherited from subquotients, so that we are finished.

For a tamely ramified niveau 2 representation the proof is essentially identical—one of the weights predicted in Conjecture 3.1 embeds in the module corresponding to the weight predicted by Serre's conjecture.

If $\rho$ is wildly ramified then we have

$$\rho = \begin{pmatrix} \omega^\alpha & * \\ & \omega^\beta \end{pmatrix}.$$

Conjecture 3.1 then predicts a weight of $F(\alpha-1,\beta)'$, which is unambiguously defined as long as $\alpha \not\equiv \beta+1 \pmod{p-1}$.

In order to apply Serre's conjecture we normalize so that $1 \le \alpha \le p-1$, and $0 \le \beta \le p-2$.

If $\alpha > \beta$ and $\alpha \neq \beta+1$, then Serre's conjecture predicts a weight of $F(\alpha-1+\beta p, 0)$, which contains $F(\alpha-1,\beta)$ as a submodule, hence we are finished as before.

If $\alpha \le \beta$, then Serre's conjecture predicts a weight of $1+\beta+p\alpha$, and we have that $F(\alpha-1,\beta)' = F(\alpha-1+p-1,\beta)$. Using [10, Table 1] as before, we find that $F(\alpha-1+p-1,\beta)$ is a subquotient of $F(\beta-1+p\alpha,0)$, which is the module corresponding to the weight predicted by Serre's conjecture. Hence, by Lemma 3.9, we are finished. ☐

**Lemma 3.9.** *If $\alpha$ is an eigenclass in $H^1(\Gamma_0(N), A)$, where $A$ is a subquotient of a $GL_2(F)$-module $B$, $\alpha$ is not a twist of a punctual eigenclass, and $p > 3$, then there is an eigenclass in $H^1(\Gamma_0(N), B)$ with the same eigenvalues as $\alpha$.*

*Proof.* Let $T \subset S \subset B$, with $S/T \cong A$, and examine the long exact cohomology sequence arising from the short exact sequence

$$0 \to T \to S \to A \to 0.$$

Note that since $p > 3$, $H^2(\Gamma_0(N), T) = 0$, so that the eigenclass $\alpha$ must come from a class $\sigma$ in $H^1(\Gamma_0(N), S)$. By [5], we may replace $\sigma$ by an eigenclass having the same eigenvalues as $\alpha$, (calling the new class $\sigma$ again). The long exact cohomology sequence arising from the short exact sequence

$$0 \to S \to B \to B/S \to 0$$

then shows that $\sigma$ goes to a nonzero class $\beta$ in $H^1(\Gamma_0(N), B)$, since it cannot come from $H^0(\Gamma_0(N), B/S)$ (as it is not a twist of a punctual class). Clearly $\beta$ has the same eigenvalues as $\sigma$. ☐

**Theorem 3.10.** *Assume that $\rho : G_{\mathbb{Q}} \to GL_n(\bar{\mathbb{F}}_p)$ is attached to an eigenclass $\alpha$ in $H^i(\Gamma_0(N), V(\epsilon))$, where $N$, $\epsilon$, and $V$ are the level, nebentype, and a weight predicted for $\rho$. Then $\rho^\vee = {}^t\rho^{-1}$ is attached to a cohomology class $\beta$ in $H^i(\Gamma_0(N), W(\epsilon^{-1}))$, where $W = V^* \otimes \det^{-(n-1)}$. Further, the level, nebentype, and a weight predicted for $\rho^\vee$ are $N$, $\epsilon^{-1}$, and $W$.*

*Proof.* The proof that there is a $\beta$ in the indicated cohomology group with $\rho^\vee$ attached is exactly the same as the proof of [4, Proposition 2.8]. For $\rho$ of niveau 1, Ash and Sinnott also prove that the invariants of $\rho^\vee$ are as above. The level and nebentype computations remain the same regardless of the niveau of the representation, so we need only show that $W$ is a predicted weight for $\rho^\vee$.

We will show that if $(b_1, \ldots, b_n)$ is a derived $n$-tuple for $\rho$, then $(-b_n, \ldots, -b_1)$ is a derived $n$-tuple for $\rho^\vee$. Then, since $(F(\alpha_1, ..., \alpha_n)')^* = F(-\alpha_n, \ldots, -\alpha_1)'$, it follows that if $V$ is a predicted weight for $\rho$, then $W$ is a predicted weight for $\rho^\vee$.

It is an easy exercise to reduce the question to simple representations of $G_p$. Suppose that $\rho$ is a simple representation of $G_p$, with the $n$-tuple $(b_1, \ldots, b_n)$ derived from it. Then there must be some exponent $m$ such that

$$\rho|_{I_p} = \begin{pmatrix} \varphi_1 & & \\ & \ddots & \\ & & \varphi_n \end{pmatrix},$$

where $(\varphi_1, \ldots, \varphi_n)$ is some permutation of $\psi_{n,1}^m, \ldots, \psi_{n,n}^m$. Then $-m$ is an exponent associated to $\rho^\vee$ in the same way, as is any multiple of $-m$ by a power of $p$. Now $m = a_1 + a_2 p + \cdots + a_n p^{n-1}$, where the $a_i$ are some permutation of the decreasing $n$-tuple $(b_1, \ldots, b_n)$, with $0 \leq a_i - a_n \leq p-1$. Let $a_k$ be the largest of the $a_i$, which is equal to $b_1$. Then $-p^{n-1-k} m$ is congruent (modulo $(p^n - 1)$) to

$$-a_{k+1} - \ldots - a_n p^{n-2-k} - a_1 p^{n-1-k} - \ldots - a_k p^{n-1},$$

with $0 \leq a_i - a_k \leq p - 1$, so that $(-b_n, \ldots, -b_1)$ is easily seen to be an $n$-tuple associated with $\rho^\vee$. □

3.1. **Heuristic for the niveau $n$ case.** For the most part, we have derived our conjecture using Serre's conjecture as a model. We can provide a suggestive heuristic for one feature of our conjecture: the weight of a niveau $n$ representation into $\mathrm{GL}_n(\bar{\mathbb{F}}_p)$.

Let $\rho : G_{\mathbb{Q}} \to \mathrm{GL}_n(\bar{\mathbb{F}}_p)$ be given such that

$$\rho|_{I_p} \sim \begin{pmatrix} \varphi_1 & & \\ & \ddots & \\ & & \varphi_n \end{pmatrix},$$

where the $\varphi_i$ are powers of a fundamental character of niveau $n$, and are conjugate to each other.

Let us suppose that $\rho$ lifts to a $p$-adic representation $\Theta$ unramified at almost all primes. Further, suppose that $\Theta$ comes from a motive $M$ with good reduction at $p$, which would conjecturally be the case $\Theta$ were attached to an automorphic representation $\pi$ of cohomological type of level $N$ prime to $p$ cf. [9]. Then $\Theta$ is crystalline. So by analogy it is reasonable to assume that $\rho$ is "crystalline" in the sense of [14], i.e. that it corresponds to a filtered Frobenius module for $\mathbb{F}_p$.

Now write $\varphi_1 = \psi^{a_1 + a_2 p + \cdots + a_n p^{n-1}}$ with $0 \leq a_i \leq p - 1$. By [14, Theorem 0.8], there is indeed a unique filtered Frobenius module $\Phi$ over $\mathbb{F}_p$ that corresponds to a representation of $G_{\mathbb{Q}_p}$ into $\mathrm{GL}_n(\mathbb{F}_p)$ whose restriction to $I_p$ is equivalent to $\rho|_{I_p}$. This is our motivation for choosing the $a_i$ in the given range.

Assuming again that $\Theta$ and $M$ exist, the Hodge numbers of $M$ would be the same as the Hodge-Tate numbers of $\Theta|_{G_{\mathbb{Q}_p}}$, and these in turn would be the same as the jumps in the filtration of the filtered Frobenius module associated to $\Theta|_{G_{\mathbb{Q}_p}}$. If we take the latter to be the same as the jumps in $\Phi$, they are $a_1, \ldots, a_n$.

Now suppose we are in the "generic case," so $|a_i - a_j| > 1$ for $i \neq j$. Let $\{b_1, \ldots, b_n\} = \{a_1, \ldots, a_n\}$, with $b_1 > b_2 > \ldots > b_n$. Assuming the general picture of Clozel (following Langlands) of the relationship between automorphic representations and motives, as found in [9], especially chapters 3 and 4, the motive $M$ predicts the existence of a $\pi$ attached to $M$ such that $\pi_\infty \otimes W$ has $(\mathfrak{g}, K)$ cohomology, where $W$ is the irreducible representation of $\mathrm{GL}_n(\mathbb{C})$ with highest weight $(b_1 - (n-1), b_2 - (n-2), \dots, b_n)$.

By analogy, we conjecture that $\rho$ will be attached to a cohomology class with weight $V = F(b_1 - (n-1), b_2 - (n-2), \dots, b_n)$. After all, $\rho$ is the reduction of $\Theta$ modulo $p$, and $W$ mod $p$ (or more precisely the reduction modulo $p$ of a model for $W$ over $\mathbb{Z}_p$) has $V$ as a composition factor. If we now require our conjecture to be closed under twisting by powers of $\omega$, a simple exercise yields the conjectural weights in Conjecture 3.1 for niveau $n$, dimension $n$, in the generic case. By "continuity" we extend the heuristic to the non-generic case.

## 4. Symmetric Squares

Using work of Ash and Tiep [7], who proved that certain Galois representations are in fact attached to cohomology eigenclasses, we are able to verify certain special cases of Conjecture 3.1.

**Theorem 4.1.** *Let $\sigma : G_{\mathbb{Q}} \to GL_2(\mathbb{F}_p)$ be a continuous irreducible odd Galois representation ramified only at p, for which Serre's conjecture is true, and let k be the weight predicted by Serre's conjecture. Then if $2 < k < (p+3)/2$, $\mathrm{Sym}^2 \sigma$ is attached to a cohomology eigenclass in weight $F(2(k-2), k-2, 0)$, and this weight is predicted by Conjecture 3.1.*

*Proof.* By [6, Proposition 2.5], we see that $\sigma$ is attached to a cohomology eigenclass in $H^1(\mathrm{SL}_2(\mathbb{Z}), U_h(\bar{\mathbb{F}}_p))$, where $h = k-2$, and $U_h(\bar{\mathbb{F}}_p) = \mathrm{Sym}^h(\bar{\mathbb{F}}_p^2)$, with the standard action of $GL_2(\bar{\mathbb{F}}_p)$ on $\bar{\mathbb{F}}_p^2$. Then, by [7, Corollary 5.3], $\mathrm{Sym}^2 \sigma$ is attached to a cohomology eigenclass in $H^3(\mathrm{SL}_3(\mathbb{Z}), F(2h, h, 0))$. Hence we need only show that the weight $F(2h, h, 0)$ is predicted by Conjecture 3.1.

If $\sigma$ has niveau 1, this is trivial, since we must have

$$\sigma|_{I_p} \sim \begin{pmatrix} \omega^{k-1} & * \\ & 1 \end{pmatrix},$$

so that

$$\mathrm{Sym}^2 \sigma|_{I_p} \sim \begin{pmatrix} \omega^{2(k-1)} & * & * \\ & \omega^{k-1} & * \\ & & 1 \end{pmatrix}.$$

If $\sigma$ has niveau 2, then we must have

$$\sigma|_{I_p} \sim \begin{pmatrix} \psi^{k-1} & \\ & \psi'^{k-1} \end{pmatrix},$$

with $1 \leq k-1 \leq (p-1)/2$, so that

$$\mathrm{Sym}^2 \sigma|_{I_p} \sim \begin{pmatrix} \psi^{2(k-1)} & & \\ & \omega^{k-1} & \\ & & \psi'^{2(k-1)} \end{pmatrix},$$

with $2 \leq 2(k-1) \leq p-1$. Clearly a predicted weight for this representation is $F(2(k-2), k-2, 0)$.

□

*Example* 4.2. Let $K$ be a totally complex $S_4$-extension of $\mathbb{Q}$, such that the quartic subfield of $K$ has discriminant $p^3$, where $p$ is a prime congruent to 5 mod 8 (for examples of such fields, see [11]. The unique three-dimensional irreducible unimodular mod $p$ representation of $S_4$ gives rise to an irreducible unimodular representation $\rho : G_{\mathbb{Q}} \to \mathrm{GL}_3(\mathbb{F}_p)$ which is ramified only at $p$. This representation is (up to a twist by a power of the cyclotomic character) the symmetric square of a two-dimensional irreducible representation $\sigma : G_{\mathbb{Q}} \to \mathrm{GL}_2(\bar{\mathbb{F}}_p)$ with projective image isomorphic to $S_4$ and image of order 96 [19]. Serre's conjecture is true for $\sigma$, since $\sigma$ has a lift to a two-dimensional irreducible complex Galois representation with solvable image, to which we may apply Langlands-Tunnell. Hence $\sigma$ is modular and so, by the $\epsilon$-conjecture, Serre's conjecture holds for $\sigma$ (see [11] for more details). One easily checks that $\sigma$ has niveau 1, and that the weight predicted by Serre's conjecture for $\sigma$ is $(p+3)/4$, so that Theorem 4.1 applies. Hence, at least one of the weights predicted for by Conjecture 3.1 yields an eigenclass with $\rho$ attached. In fact, this weight is $F((p-5)/2,(p-5)/4,0) \otimes \det^{3(p-1)/4}$.

*Example* 4.3. Let $K$ be a totally complex $S_4$-extension of $\mathbb{Q}$, such that the quartic subfield of $K$ has discriminant $-p$, where $p$ is a prime congruent to 3 mod 8. Let $\rho$ be the unimodular irreducible three-dimensional Galois representation associated to $K$ as above. Again, there is a two-dimensional irreducible representation $\sigma : G_{\mathbb{Q}} \to \mathrm{GL}_2(\mathbb{F}_p)$ with projective image isomorphic to $S_4$, such that $\sigma$ is ramified only at $p$ [19] (this time the image of $\sigma$ has order 48), and (again up to a twist by a power of the cyclotomic character), $\rho$ is the symmetric square of $\sigma$ (note that up to twisting the symmetric square depends only on the projectivization of a representation). One checks easily that Serre's conjecture predicts a weight of $(p+1)/2$ for $\sigma$ (again $\sigma$ has niveau 1), and that (just as above), Serre's conjecture is true for $\sigma$. Hence, one of the weights predicted by Conjecture 3.1 does in fact contain an eigenclass with $\rho$ attached. In this case, the weight is $F(p-3,(p-3)/2,0) \otimes \det^{(p-1)/2}$.

*Example* 4.4. Let $K$ be a complex $S_4$ extension of $\mathbb{Q}$, with $K$ ramified at only one prime $p$, with $p$ congruent to 3 modulo 8, and with ramification index at $p$ equal to 4 (for examples of such extensions, see [11]). Let $\rho$ be the unique unimodular irreducible three-dimensional mod $p$ Galois representation with image isomorphic to $S_4$ and such that the fixed field of the kernel of $\rho$ is $K$. Then, up to twisting, $\rho$ is the symmetric square of a representation $\sigma : G_{\mathbb{Q}} \to \mathrm{GL}_2(\mathbb{F}_p)$ with image isomorphic to $\tilde{S}_4$. In this case, $\sigma$ has niveau 2, Serre's conjecture is true for $\sigma$ and its twists, and a twist of $\sigma$ has weight $(p+5)/4$ [11], so that Theorem 4.1 applies. Hence, one of the weights predicted for $\rho$ gives a cohomology group which contains an eigenclass predicted for $\rho$. In this case, the weight which works is $F((p-3)/2,(p-3)/4,0) \otimes \det^{(3p-5)/4}$.

## 5. Niveau 1 representations

5.1. **Reducible representations in level 1.** In [4] Ash and Sinnott deal extensively with reducible representations ramified at only one prime. Each of their examples was a direct sum of an even two-dimensional representation with a one-dimensional representation, and they included cases where the two-dimensional representation had image isomorphic to a dihedral group, or projective image isomorphic to $A_4$. They did not give examples in which the projective image was isomorphic to $S_4$ or $A_5$.

We recall their construction from [4]:

Let $\sigma : G_{\mathbb{Q}} \to GL_2(\bar{\mathbb{F}}_p)$ be an irreducible representation with the following properties:

1. $\sigma$ is unramified outside $p$;
2. The image of $\sigma$ has order relatively prime to $p$;
3. $\sigma(\mathrm{Frob}_\infty)$ is central, where $\mathrm{Frob}_\infty$ is a complex conjugation in $G_{\mathbb{Q}}$;
4. $\sigma(G_{p,0})$ has order dividing $p-1$.

Then choosing integers $j$ and $k$ appropriately, we find that the representation $\rho = \sigma \otimes \omega^j \oplus \omega^k$ is three-dimensional and odd, and by adjusting $j$ and $k$, we may adjust the predicted weight of $\rho$ to some extent. In particular, we need to choose $j$ and $k$ to have opposite parity if $\sigma(\mathrm{Frob}_\infty) = 1$, and the same parity if $\sigma(\mathrm{Frob}_\infty) = -1$. In addition, we choose $j$ and $k$ to give the simplest possible weight.

The reducibility of these representations makes it possible to reduce the weight to calculable levels, however in the examples that we consider here the weight will still be quite high. Hence, rather than being able to calculate many Hecke eigenvalues, time constraints make it impractical to calculate more eigenvalues than those at 2 and 3.

We begin by specifying the fixed field of the kernel of the projective image of $\sigma$, which will be a totally real number field.

5.1.1. *Representations of type $A_4$.* In [4] Ash and Sinnott present several examples of reducible Galois representations which are sums of one-dimensional characters with even two-dimensional representations having projective image isomorphic to $A_4$. Using the same computational techniques as in [1], we have been able to find other $A_4$-extensions for which we can compute the predicted quasicuspidal homology classes. These examples are given in Table 1.

We begin with a quartic polynomial $f$ which has four real roots and whose splitting field $K$ is an $A_4$-extension of $\mathbb{Q}$ ramified only at one prime $p$. We know (by the lifting Lemma in [4]) that $K$ sits inside an $\hat{A}_4$-extension $\hat{K}$ of $\mathbb{Q}$, with $\hat{K}/\mathbb{Q}$ ramified only at $p$. In fact, there are two possibilities for $\hat{K}$, following [4], we take $\hat{K}$ to be the one which has ramification index 3 at $p$. Let $K_4$ be the quartic extension of $\mathbb{Q}$ defined by $f$. We note that $K_4$ must be contained in an octic subextension $K_8$ of $\hat{K}$, with $K_8/K_4$ unramified at all finite primes. Since $K_8$ has $\hat{K}$ as its Galois closure, we may determine whether $\hat{K}$ is totally real or totally complex by comparing the two-ranks of the class group and narrow class group of $K_4$. For instance, when $p = 1009$, the class number of $K_4$ is two, and the narrow class group is cyclic of order four. Thus, the two class groups have the same two-rank, so $\hat{K}$ must be real (since $K_8$ and all its conjugates must be real). If $\hat{K}$ is totally real, we write its sign as 1, otherwise its sign is $-1$.

Now $\hat{A}_4$ has a unique two-dimensional irreducible unimodular mod $p$ representation $\sigma : G_{\mathbb{Q}} \to GL_2(\mathbb{F}_p)$. We see easily that $\sigma|_{I_p} = \omega^d \oplus \omega^{-d}$, with $d = (p-1)/3$. We now take $\rho = \sigma \otimes \omega^j \oplus \omega^k$, with $j = 2d$ and $k = 1$ if the sign of $\hat{K}$ is 1, and $k = 2$ otherwise. We note that $\rho$ satisfies the conditions of the construction of [4], and has a predicted weight of $F((p-1)/3-2, 0, 0)$ if $k = 1$, and $F((p-1)/3-2, 1, 0)$ if $k = 2$.

For each of the examples in Table 1 we have calculated the interior homology of $SL_3(\mathbb{Z})$ in the given weight, using the techniques described in [1], and found it to be one-dimensional. We have also calculated the Hecke eigenvalues at 2 and 3, and found that they exactly match the values predicted from the characteristic polynomial of the image of Frobenius under $\rho$ by Conjecture 3.1.

| Polynomial | sign | $p$ | k | Weight |
|---|---|---|---|---|
| $x^4-2x^3-13x^2-9x+4$ | -1 | 163 | 2 | $F(52,1,0)$ |
| $x^4-x^3-16x^2+3x+1$ | 1 | 277 | 1 | $F(90,0,0)$ |
| $x^4-x^3-10x^2+3x+20$ | -1 | 349 | 2 | $F(114,1,0)$ |
| $x^4-x^3-13x^2+12x+16$ | -1 | 397 | 2 | $F(130,1,0)$ |
| $x^4-2x^3-19x^2+29x+1$ | -1 | 547 | 2 | $F(180,1,0)$ |
| $x^4-2x^3-31x^2-51x-4$ | 1 | 607 | 1 | $F(200,0,0)$ |
| $x^4-2x^3-39x^2+x+125$ | 1 | 1009 | 1 | $F(334,0,0)$ |
| $x^4-2x^3-51x^2+100x+83$ | 1 | 1399 | 1 | $F(464,0,0)$ |
| $x^4-2x^3-51x^2+32x+192$ | 1 | 1699 | 1 | $F(564,0,0)$ |
| $x^4-2x^3-37x^2+10x+29$ | 1 | 1777 | 1 | $F(590,0,0)$ |
| $x^4-2x^3-43x^2+127x-55$ | 1 | 1951 | 1 | $F(648,0,0)$ |

TABLE 1. Reducible representations of type $A_4$

5.1.2. *Representations of type $S_4$.* Totally real $S_4$-extensions ramified at only one prime can have two types of ramification; either the ramification index is 2 or the ramification index is 4. For our purposes, the extensions with ramification index 4 are better (since they yield lower weights), although they are more difficult to find. They can, however, be found by application of explicit class field theory, and many such examples are known. Only the two below yield predicted weights which are feasible for computation.

*Example* 5.1. Let $K$ be the splitting field of the polynomial $x^4 - x^3 - 1017x^2 + 9665x+60608$. Then $K$ is a totally real $S_4$-extension of $\mathbb{Q}$, ramified only at $p = 2713$, with ramification index $e = 4$. Let $\tilde{S}_4$ be the central extension of $S_4$ by $\mathbb{Z}/2\mathbb{Z}$ which is isomorphic to $\mathrm{GL}_2(\mathbb{F}_3)$. Then $K$ embeds in an $\tilde{S}_4$-extension $\tilde{K}$ of $\mathbb{Q}$ (by the lifting lemma of [4]), and $\tilde{K}/K$ must further ramify at $p$ (as described in [11]), so that in $\tilde{K}$, $p = 2713$ has $e = 8$. We need to determine whether $\tilde{K}$ is totally real or totally complex. To do this, we note that $\tilde{S}_4$ has three conjugacy classes of subgroups of order 6, and that each subgroup of order 12 contains exactly one subgroup of order 6 from each conjugacy class. In terms of field extensions then, each subfield of $\tilde{K}$ of degree 4 has exactly three quadratic extensions lying in $\tilde{K}$. Hence, if $\tilde{K}$ is totally real, the degree 4 subfield $K_4$ of $K$ must have a Klein four extension contained inside $\tilde{K}$, hence ramified only at $p$ (in particular not ramified at infinity). Such an extension would lie inside the ray class field of $K_4$ modulo $\mathfrak{p}^m$ (where $\mathfrak{p}$ is the unique prime of $K_4$ lying over $p$). However, the two part of the ray class group of $K_4$ modulo $\mathfrak{p}^m$ is cyclic for every $m$ [17]. Hence, $\tilde{K}$ must be totally complex.

Now, we let $\sigma$ be the two-dimensional representation $\sigma : G_{\mathbb{Q}} \to GL_2(\bar{\mathbb{F}}_p)$, with image isomorphic to $\tilde{S}_4$ and kernel equal to $G_{\tilde{K}}$, chosen such that

$$\sigma|_{I_p} = \begin{pmatrix} \omega^{3(p-1)/8} & \\ & \omega^{(p-1)/8} \end{pmatrix} = \begin{pmatrix} \omega^{3(339)} & \\ & \omega^{339} \end{pmatrix}.$$

Taking $j = -339$, $k = 1$, (with the same parity, since $\sigma(\mathrm{Frob}_\infty) = -1$), we see that $\rho = \sigma \otimes \omega^j \oplus \omega^k$ has

$$\rho|_{I_p} \sim_L \begin{pmatrix} \omega^{678} & & \\ & \omega^1 & \\ & & \omega^0 \end{pmatrix},$$

where we take $L = \begin{pmatrix} * & & * \\ & * & \\ * & & * \end{pmatrix}$. Then the weight predicted by Conjecture 3.1 is $F(678-2, 1-1, 0)' = F(676, 0, 0)$, the level is 1, and the nebentype is trivial. Computations using the techniques of [1] show that the interior cohomology is in fact one-dimensional. The Hecke eigenvalues at 2 and 3 correspond exactly to $\sigma$, as predicted by Conjecture 3.1.

A similar construction can be performed with the splitting field $K$ of the polynomial $x^4 - 6668x^3 + 16598046x^2 - 18278822428x + 7514424150025$, which is a totally real $S_4$-extension of $\mathbb{Q}$ ramified only at $p = 3137$. In this case $\tilde{K}$ is totally real, and the predicted weight is $F(782, 0, 0)$. Again, the homology is one-dimensional, and the eigenvalues at 2 match $\rho$. The image of the Frobenius at 3 is of order 8, however, and presents some difficulty. We have determined $\sigma$ (and hence also $\rho$) by a local condition at $p$, namely its restriction to inertia at $I_p$. Determining the Frobenius at 3 is a local condition at 3, and combining these two determinations (in order to determine $\mathrm{Tr}(\rho(\mathrm{Frob}_3))$) is a global problem which involves calculations in a large number field. We thus have two possibilities for eigenvalues at 3 which would correspond to $\rho$, which we cannot distinguish. One of these possibilities does in fact occur in the predicted cohomology.

5.1.3. *Representations of type $A_5$.* Ash and Sinnott's construction works best with $A_5$-extensions if the ramification index of the unique ramified prime is as large as possible. However, totally real $A_5$-extensions of $\mathbb{Q}$ ramified at only one prime with ramification index 5 are quite difficult to find. The second author thanks Steve Harding for showing him the example below with $p = 3821$.

*Example* 5.2. Let $K$ be the splitting field of the polynomial $f = x^5 - 7402x^3 - 3701x^2 + 14804x + 11103$. Then $K$ is a totally real Galois extension of $\mathbb{Q}$, with Galois group $A_5$, ramified only at $p = 3701$. $K$ must lie inside an extension $\hat{K}$ of $\mathbb{Q}$ with Galois group $\hat{A}_5$ (the unique non-split central extension of $A_5$ by $\mathbb{Z}/2\mathbb{Z}$). In fact, $K$ lies inside two such extensions, one in which primes above $p$ ramify further, and one in which primes above $p$ do not ramify further.

Let $\hat{K}$ be an $\hat{A}_5$-extension of $\mathbb{Q}$, containing $K$, in which $p$ has ramification index 5. Let $H$ be a subgroup of $\hat{A}_5$ of order 20. Using the computer algebra system Magma, one can see that $H$ has a quotient group which is cyclic of order 4. Hence, the degree 6 subextension of $K$ must have a cyclic quartic extension contained in $\hat{K}$, which is unramified at all finite primes.

A defining polynomial for the degree 6 subextension of $K$ may be found as the minimal polynomial of the element

$$\alpha_1\alpha_2 + \alpha_2\alpha_3 + \alpha_3\alpha_4 + \alpha_4\alpha_5 + \alpha_5\alpha_1,$$

where $\alpha_i$, $1 \le i \le 5$ are the roots of $f$. Using GP/PARI to compute the ideal class group and the narrow class group, we find that both are cyclic of order 4. Hence, the only possible cyclic quartic extension of the degree 6 subfield of $K$ which is unramified at all finite primes is also unramified at infinity, so that $\hat{K}$ is totally real.

Now $\hat{A}_5$ has two two-dimensional mod $p$ representations. Call them $\sigma$ and $\sigma'$. On inertia at $p = 3701$, we may choose $\sigma$ and $\sigma'$ so that

$$\sigma|_{I_p} \sim \begin{pmatrix} \omega^{3(p-1)/5} & \\ & \omega^{2(p-1)/5} \end{pmatrix}, \quad \text{and} \quad \sigma'|_{I_p} \sim \begin{pmatrix} \omega^{(p-1)/5} & \\ & \omega^{-(p-1)/5} \end{pmatrix}.$$

If we let $\rho = \sigma \otimes \omega^{-2(p-1)/5} \oplus \omega$, then $\rho$ is an odd three-dimensional representation, and if it is conjugated to land inside

$$L = \begin{pmatrix} * & & * \\ & * & \\ * & & * \end{pmatrix}$$

then it satisfies strict parity. We then have

$$\rho|_{I_p} \sim_L \begin{pmatrix} \omega^{(p-1)/5} & & \\ & \omega & \\ & & 1 \end{pmatrix},$$

giving a predicted weight of $F((p-1)/5 - 2, 0, 0) = F(738, 0, 0)$. We may calculate the quasicuspidal homology in this weight, and find that it is one-dimensional, and has the appropriate eigenvalues at 2 and 3 to correspond to $\rho$. In this case, there is an ambiguity similar to that in the preceding example, in that we cannot determine which of the two conjugacy classes of order 5 contains the Frobenius at 2. The computed eigenvalues at 2 are in fact one of the two possible pairs of values which could have $\rho$ attached.

A similar calculation may be carried out for the $A_5$-extension defined by the polynomial $x^5 - 3821x^3 - 3821x^2 + 3821x + 3821$, and ramified only at $p = 3821$. In this case $\hat{K}$ is again totally real. Hence, as above, we get a predicted weight for $\rho$ of $F((p-1)/5 - 2, 0, 0) = F(762, 0, 0)$. Calculating the quasicuspidal homology in this weight yields a one-dimensional space, which has appropriate eigenvalues at 2 and 3 to correspond to $\rho$, with the ambiguity that we are unable to determine the conjugacy class of elements of order 10 (resp. 5) containing the Frobenius at 2 (resp. 3), just as in the previous examples.

5.2. **Reducible representations in higher level.** With the introduction of levels higher than one, we gain immensely in reducing the weight of the representations that we can find. In particular, we find that we can actually compute "companion forms," or classes with different weights, attached to the same representation. These offer important examples of Conjecture 3.1.

We work out one interesting example in full detail, and describe others in a table format.

*Example* 5.3. Let $K$ be the $S_3$-extension of $\mathbb{Q}$ given as the splitting field of the polynomial $x^3 - x^2 - 3x + 1$. Then $K$ is ramified only at $p = 37$ (with ramification index 2) and at $q = 2$ (with ramification index 3). Since $S_3$ has a two-dimensional mod 37 representation, we obtain a two-dimensional Galois representation $\sigma : G_{\mathbb{Q}} \to \mathrm{GL}_2(\mathbb{F}_{37})$, with the fixed field of the kernel of $\sigma$ equal to $K$. Let $\omega$ be the cyclotomic character modulo 37, and let $\rho = \sigma \oplus \omega$. We note that $\sigma$ is an even representation, since $K$ is totally real, so we want to conjugate $\rho$ to land inside the Levi subgroup

$$L = \begin{pmatrix} * & & * \\ & * & \\ * & & * \end{pmatrix}.$$

Now inside $L$, the image of complex conjugation is conjugate to the matrix

$$\begin{pmatrix} 1 & & \\ & -1 & \\ & & 1 \end{pmatrix},$$

so that $\rho$ satisfies strict parity.

One sees easily that the level of $\rho$ is equal to the level of $\sigma$, which is $2^2$ (since the ramification at 2 is tame and the image of inertia at 2 under $\sigma$ does not fix a subspace). The nebentype of $\rho$ is trivial, since the determinant is just $\omega^{19}$. Finally, if we examine the restriction of $\rho$ to inertia at 37, we find that,

$$\rho|_{I_{37}} \sim_L \begin{pmatrix} \omega^{18} & & \\ & \omega^1 & \\ & & \omega^0 \end{pmatrix}.$$

Thus, the weight predicted by Conjecture 3.1 is $F(18-2, 1-1, 0)' = F(16, 0, 0)$. When we compute the cohomology in this weight, in level 4 with trivial nebentype, we obtain a fifteen-dimensional space, containing a one-dimensional eigenspace with eigenvalues given by the following table:

| Eigenvalues | 2 | 3 | 5 | 7 | 11 | 13 | 17 | 19 | 23 | 29 | 31 | 37 | 41 | 43 | 47 |
|---|---|---|---|---|---|---|---|---|---|---|---|---|---|---|---|
| $a(\ell, 1)$ | * | 2 | 5 | 6 | 10 | 13 | 17 | 19 | 23 | 29 | 31 | 1 | 3 | 6 | 9 |
| $a(\ell, 2)$ | * | 24 | 22 | 15 | 26 | 17 | 13 | 35 | 8 | 29 | 31 | 1 | 27 | 6 | 25 |

We now compute the trace $\mathrm{Tr}(\rho(\mathrm{Frob}_\ell))$ and $T_2(\rho(\mathrm{Frob}_\ell))$ (the sum of products of pairs of eigenvalues) for $\ell$ between 2 and 47. To do this we note that the characteristic polynomial of $\rho$ is

$$\begin{aligned} \det(I - x\rho(\mathrm{Frob}_\ell)) &= \det(I - x\sigma(\mathrm{Frob}_\ell))(1 - x\omega(\mathrm{Frob}_\ell)) \\ &= (1 - \mathrm{Tr}(\sigma(\mathrm{Frob}_\ell))x + D(\sigma(\mathrm{Frob}_\ell))x^2)(1 - \ell x) \\ &= 1 - (\mathrm{Tr}(\sigma(\mathrm{Frob}_\ell)) + \ell)x + \\ &\qquad (D(\sigma(\mathrm{Frob}_\ell)) + \mathrm{Tr}(\sigma(\mathrm{Frob}_\ell))\ell))x^2 - D(\sigma(\mathrm{Frob}_\ell))\ell x^3 \end{aligned}$$

so that the trace of $\rho(\mathrm{Frob}_\ell)$ is $\mathrm{Tr}(\sigma(\mathrm{Frob}_\ell)) + \ell$ and $T_2(\rho(\mathrm{Frob}_\ell)) = D(\sigma(\mathrm{Frob}_\ell)) + \mathrm{Tr}(\sigma(\mathrm{Frob}_\ell)\ell)$. Using GP/PARI, we may calculate these two values for $\ell$ from 2 to 47, (excluding the ramified primes 2 and 37), and we find that they exactly match the values of $a(\ell, 1)$ and $\ell a(\ell, 2)$, calculated above.

Other reducible examples are easily computed just as above. In each row of Table 2 we give a polynomial whose splitting field $K$ is a totally real Galois extension of $\mathbb{Q}$ with Galois group $G$, such that $G$ has a unique two-dimensional representation $\sigma$ modulo $p$. We also give the predicted weight(s), level, and nebentype of the cohomology class corresponding to $\rho = \sigma \oplus \omega$. Several examples have more than one predicted weight, coming from multiple orderings of the diagonal characters. Such predictions actually occur in all of these examples, but most are too large for us to calculate. For all of the examples in this table, all Hecke eigenvalues up to $\ell = 47$ coincided exactly with the coefficients of the characteristic polynomial of the image of Frobenius, as predicted by Conjecture 3.1.

| Polynomial | $G$ | $p$ | Weight(s) | Level | $\epsilon$ |
|---|---|---|---|---|---|
| $x^3 - x^2 - 3x + 1$ | $S_3$ | 37 | $F(16, 0, 0)$ | 4 | 1 |
| $x^3 - x^2 - 4x + 2$ | $S_3$ | 79 | $F(37, 0, 0)$ | 4 | $\epsilon_4$ |
| $x^3 - x^2 - 5x - 1$ | $S_3$ | 101 | $F(48, 0, 0)$ | 4 | 1 |
| $x^3 - x^2 - 4x + 1$ | $S_3$ | 107 | $F(51, 0, 0)$ | 3 | $\epsilon_3$ |
| $x^3 - x^2 - 5x + 4$ | $S_3$ | 67 | $F(31, 0, 0)$ | 7 | $\epsilon_7$ |
| $x^3 - 5x - 1$ | $S_3$ | 43 | $F(19, 0, 0)$ | 11 | $\epsilon_{11}$ |
| | $S_3$ | 11 | $F(3, 0, 0)$ | 43 | $\epsilon_{43}$ |
| $x^3 - 7x - 5$ | $S_3$ | 41 | $F(18, 0, 0)$ | 17 | $\epsilon_{17}$ |
| | | 17 | $F(6, 0, 0)$ | 41 | $\epsilon_{41}$ |
| $x^3 - x^2 - 6x + 5$ | $S_3$ | 5 | $F(0, 0, 0)$ | 157 | $\epsilon_{157}$ |
| $x^3 - 7x - 1$ | $S_3$ | 5 | $F(0, 0, 0)$ | 269 | $\epsilon_{269}$ |
| $x^3 - x^2 - 9x + 8$ | $S_3$ | 7 | $F(8, 6, 2)$, $F(6, 6, 4)$ | 53 | $\epsilon_{53}$ |
| $x^4 - x^3 - 3x^2 + x + 1$ | $D_4$ | 5 | $F(0, 0, 0)$, $F(6, 4, 2)$ | 29 | $\epsilon_{29}$ |
| $x^4 - x^3 - 3x^2 + x + 1$ | $D_4$ | 29 | $F(12, 0, 0)$ | 5 | $\epsilon_5$ |
| $x^4 - x^3 - 5x^2 + 2x + 4$ | $D_4$ | 5 | $F(0, 0, 0)$, $F(6, 4, 2)$ | 89 | $\epsilon_{89}$ |
| | | 89 | $F(42, 0, 0)$ | 5 | $\epsilon_5$ |
| $x^4 - 2x^3 - 4x^2 + 5x + 5$ | $D_4$ | 5 | $F(0, 0, 0)$, $F(6, 4, 2)$ | 101 | $\epsilon_{101}$ |
| | | 101 | $F(48, 0, 0)$ | 5 | $\epsilon_5$ |
| $x^4 - x^3 - 7x^2 + 3x + 9$ | $D_4$ | 5 | $F(0, 0, 0)$ | 181 | $\epsilon_{181}$ |
| | | 181 | $F(88, 0, 0)$ | 5 | $\epsilon_5$ |
| $x^4 - 2x^3 - 4x^2 + 5x + 2$ | $D_4$ | 17 | $F(6, 0, 0)$ | 53 | $\epsilon_{53}$ |
| | | 53 | $F(24, 0, 0)$ | 17 | $\epsilon_{17}$ |
| $x^4 - x^3 - 6x^2 + 8x - 1$ | $D_4$ | 13 | $F(4, 0, 0)$ | 61 | $\epsilon_{61}$ |
| | | 61 | $F(28, 0, 0)$ | 13 | $\epsilon_{13}$ |
| $x^4 - x^3 - 5x^2 + x + 1$ | $D_4$ | 13 | $F(4, 0, 0)$ | 53 | $\epsilon_{53}$ |
| | | 53 | $F(24, 0, 0)$ | 13 | $\epsilon_{13}$ |

TABLE 2. Reducible higher level niveau 1 examples (Even two-dimensional plus $\omega^1$)

We may also apply Conjecture 3.1 to reducible representations which are the sum of an odd two-dimensional representation and a character. In order to satisfy strict parity, such a representation must land inside a Levi subgroup of the form

$$L = \begin{pmatrix} * & & \\ & * & * \\ & * & * \end{pmatrix} \quad \text{or} \quad L = \begin{pmatrix} * & * & \\ * & * & \\ & & * \end{pmatrix}.$$

For each such three-dimensional $\rho$ we thus have four predicted weights, two from each choice of Levi subgroup. In Table 3, for each example, we give a polynomial $f$ which has Galois group $G = S_3$ or $D_4$, together with a prime $p$ and the ramification index of $p$ in the splitting field $K$ of $f$. If we let $\sigma$ be the unique two-dimensional mod $p$ Galois representation arising from $K$, and $\rho = \sigma \oplus \omega^0$, we also give the level $N$ and nebentype $\epsilon$ associated to $\rho$, and the set of predicted weights arising from Conjecture 3.1. In this case we are able to compute with all the predicted weights, and find that in every case, an eigenclass with the correct eigenvalues (up to $\ell = 47$) appears in every predicted weight.

| Galois Representation | | Weights |
|---|---|---|
| $p = 7$, $e = 3$ <br> $G = S_3$ | $N = 19$, $\epsilon = \epsilon_{19}$ <br> $x^3 - x^2 + 5x - 6$ | $F(2,1,0)$, $F(4,3,2)$ <br> $F(6,3,0)$, $F(10,7,4)$ |
| $p = 7$, $e = 3$ <br> $G = S_3$ | $N = 47$, $\epsilon = \epsilon_{47}$ <br> $x^3 - x^2 - 2x - 27$ | $F(2,1,0)$, $F(4,3,2)$ <br> $F(6,3,0)$, $F(10,7,4)$ |
| $p = 7$, $e = 3$ <br> $G = S_3$ | $N = 59$, $\epsilon = \epsilon_{59}$ <br> $x^3 - x^2 + 5x + 8$ | $F(2,1,0)$, $F(4,3,2)$ <br> $F(6,3,0)$, $F(10,7,4)$ |
| $p = 7$, $e = 3$ <br> $G = S_3$ | $N = 59$, $\epsilon = \epsilon_{59}$ <br> $x^3 - x^2 - 9x + 36$ | $F(2,1,0)$, $F(4,3,2)$ <br> $F(6,3,0)$, $F(10,7,4)$ |
| $p = 7$, $e = 3$ <br> $G = S_3$ | $N = 59$, $\epsilon = \epsilon_{59}$ <br> $x^3 - x^2 - 2x - 20$ | $F(2,1,0)$, $F(4,3,2)$ <br> $F(6,3,0)$, $F(10,7,4)$ |
| $p = 19$, $e = 3$ <br> $G = S_3$ | $N = 3$, $\epsilon = \epsilon_3$ <br> $x^3 - x^2 - 6x - 12$ | $F(10,5,0)$, $F(16,11,6)$ <br> $F(22,11,0)$, $F(34,23,12)$ |
| $p = 13$, $e = 3$ <br> $G = S_3$ | $N = 43$, $\epsilon = \epsilon_{43}$ <br> $x^3 - x^2 - 17x + 38$ | $F(6,3,0)$, $F(10,7,4)$ <br> $F(14,7,0)$, $F(22,15,8)$ |
| $p = 3$, $e = 2$ <br> $G = D_4$ | $N = 13$, $\epsilon = \epsilon_{13}$ <br> $x^4 + x^2 - 3$ | $F(2,1,1)$, $F(1,1,0)$ <br> $F(0,0,0)$ |
| $p = 3$, $e = 2$ <br> $G = D_4$ | $N = 37$, $\epsilon = \epsilon_{37}$ <br> $x^4 + 5x^2 - 3$ | $F(2,1,1)$, $F(1,1,0)$ <br> $F(0,0,0)$ |
| $p = 3$, $e = 2$ <br> $G = D_4$ | $N = 61$, $\epsilon = \epsilon_{61}$ <br> $x^4 - 7x^2 - 3$ | $F(2,1,1)$, $F(1,1,0)$ <br> $F(0,0,0)$ |
| $p = 3$, $e = 2$ <br> $G = D_4$ | $N = 73$, $\epsilon = \epsilon_{73}$ <br> $x^4 + 34x^2 - 3$ | $F(2,1,1)$, $F(1,1,0)$ <br> $F(0,0,0)$ |
| $p = 5$, $e = 2$ <br> $G = D_4$ | $N = 39$, $\epsilon = \epsilon_3\epsilon_{13}$ <br> $x^4 - x^3 - 8x - 1$ | $F(6,5,2)$ <br> $F(4,1,0)$ |

TABLE 3. Reducible higher level niveau 1 examples (Odd two-dimensional plus $\omega^0$)

The last examples in the table, in which $\sigma$ has image isomorphic to $D_4$ (the dihedral group with 8 elements) are interesting in that fewer than four weights are predicted. In these cases, the four predicted weights are not distinct, so that the total number of weights in which we expect to find eigenvalues with $\rho$ attached is less than four. For instance, in the last example in Table 3, in which $p = 5$, the image of inertia at 5 is contained in the center of $D_4$, so that the restriction of $\sigma$ to inertia at 5 has diagonal characters $\omega^2$ and $\omega^2$. The coincidence of these diagonal characters results in the fact that only two distinct weights are predicted.

5.3. **Irreducible representations in higher level.** In order to find irreducible three-dimensional Galois representations, it is necessary to find Galois groups which have irreducible three-dimensional mod $p$ representations. For $p$ larger than 3, this is easily done: the groups $A_4$, $S_4$, and $A_5$ all have three-dimensional irreducible mod $p$ representations. We will thus concentrate primarily on representations (up to a twist) whose images are isomorphic to one of these groups. Of course we deal only with odd representations. For all the irreducible niveau 1 representations presented in this section, the three-dimensional Galois representation is a symmetric square of an odd two-dimensional representation, hence the correspondences presented here are not native three-dimensional phenomena.

5.3.1. *Representations of type $A_4$.* Suppose that $p$ is a prime congruent to 1 mod 3, and that $K$ is a totally complex $A_4$-extension ramified at $p$, with ramification index 3. There may be other ramified primes, which would then contribute to the level. Since $A_4$ has an irreducible 3-dimensional mod $p$ representation, we obtain an irreducible three-dimensional representation $\rho : G_{\mathbb{Q}} \to \mathrm{GL}_3(\mathbb{F}_p)$. We observe that the restriction of $\rho$ to inertia at $p$ is:

$$\rho|_{I_p} = \begin{pmatrix} \omega^a & & \\ & \omega^b & \\ & & \omega^c \end{pmatrix},$$

where $(a, b, c)$ is some permutation of $(2(p-1)/3, (p-1)/3, 0)$. The six permutations of $(2(p-1)/3, (p-1)/3, 0)$ then give six predicted weights for $\rho$. The six weights are

| Triple | Weight |
|---|---|
| $(2(p-1)/3, (p-1)/3, 0)$ | $F(2(p-4)/3, (p-4)/3, 0)$ |
| $((p-1)/3, 0, 2(p-1)/3)$ | $F(2(p-4)/3, (p-4)/3, 0) \otimes \det^{2(p-1)/3}$ |
| $(0, 2(p-1)/3, (p-1)/3)$ | $F(2(p-4)/3, (p-4)/3, 0) \otimes \det^{(p-1)/3}$ |
| $((p-1)/3, 2(p-1)/3, 0)$ | $F(2(2p-5)/3, (2p-5)/3, 0)$ |
| $(2(p-1)/3, 0, (p-1)/3)$ | $F(2(2p-5)/3, (2p-5)/3, 0) \otimes \det^{(p-1)/3}$ |
| $(0, (p-1)/3, 2(p-1)/3)$ | $F(2(2p-5)/3, (2p-5)/3, 0) \otimes \det^{2(p-1)/3}$, |

Hence, we expect to find three cohomology eigenclasses, each with one of $\rho$, $\rho \otimes \omega^{(p-1)/3}$, and $\rho \otimes \omega^{2(p-1)/3}$ attached in each of the two weights $F(2(p-4)/3, (p-4)/3, 0)$ and $F(2(2p-5)/3, (2p-5)/3, 0)$. In fact, however, since $\rho \otimes \omega^{(p-1)/3} \sim \rho$, the three eigenclasses may coincide, and there may actually be only one such eigenclass in each weight. In practice, in order to compute the cohomology associated to a representation as above, we will often have to twist by a character which is unramified at $p$ in order to reduce the level. We illustrate with an example.

*Example* 5.4. Let $K$ be the splitting field of the polynomial $x^4 - x^3 + 5x^2 - 4x + 3$, which is ramified at 7 (with $e = 3$), and at 13 (with $e = 2$). The predicted weights are $F(2, 1, 0)$ and $F(6, 3, 0)$. The level of $\rho$ is $13^2$, and the nebentype is trivial. Unfortunately this level is too large for us to use in computations. However, $\rho \otimes \epsilon_{13}$ is easily seen to have level 13 and nebentype $\epsilon_{13}$. Thus, we predict the existence of cohomology eigenclasses in weights $F(2, 1, 0)$ and $F(6, 3, 0)$, level 13, and nebentype $\epsilon_{13}$, which are attached to $\rho \otimes \epsilon$. Direct computation shows that these eigenclasses do in fact exist, and that the eigenvalues match, at least up to $\ell = 47$.

Other $A_4$-extensions which give rise to computable cohomology classes are shown in Table 4. Each example in this table gives a polynomial $f$ with Galois group $G$. The prime $p$, together with its ramification index $e$ in the splitting field of $f$ is given. When $G$ equals $A_4$, $\rho$ is the twist by the character $\chi$ of the unique irreducible three-dimensional mod $p$ representation of $G_{\mathbb{Q}}$ cutting out the splitting field of $f$. The level $N$, nebentype $\epsilon$, and predicted weights for $\rho$ are indicated in the table. In all cases, we have computationally verified the existence of an eigenclass in the predicted weight, level, and character, with the correct eigenvalues (up to $\ell = 47$) to have $\rho$ attached.

| Galois Representation | | Predicted Weights |
|---|---|---|
| $p=7$, $e=3$<br>$G=A_4,\chi=\epsilon_{13}$ | $N=13$ $\epsilon=\epsilon_{13}$<br>$x^4-x^3+5x^2-4x+3$ | $F(2,1,0)$, $F(4,3,2)$, $F(6,5,4)$<br>$F(6,3,0)$, $F(8,5,2)$, $F(10,7,4)$ |
| $p=7$ $e=3$<br>$G=A_4,\chi=\epsilon_{29}$ | $N=29$ $\epsilon=\epsilon_{29}$<br>$x^4-x^3+5x^2-6x+7$ | $F(2,1,0)$, $F(4,3,2)$, $F(6,5,4)$<br>$F(6,3,0)$, $F(8,5,2)$, $F(10,7,4)$ |
| $p=7$ $e=3$<br>$G=A_4,\chi=1$ | $N=2^6$ $\epsilon=1$<br>$x^4-2x^3+2x^2+2$ | $F(2,1,0)$, $F(4,3,2)$, $F(6,5,4)$<br>$F(6,3,0)$, $F(8,5,2)$, $F(10,7,4)$ |
| $p=13$, $e=3$<br>$G=A_4,\chi=\epsilon_5$ | $N=5$, $\epsilon=\epsilon_5$<br>$x^4-x^3-3x+4$ | $F(6,3,0)$, $F(10,7,4)$, $F(14,11,8)$<br>$F(14,7,0)$, $F(18,11,4)$, $F(22,15,8)$ |
| $p=13$, $e=3$<br>$G=A_4,\chi=1$ | $N=5^2$, $\epsilon=1$<br>$x^4-x^3-3x+4$ | $F(6,3,0)$, $F(10,7,4)$, $F(14,11,8)$<br>$F(14,7,0)$, $F(18,11,4)$, $F(22,15,8)$ |
| $p=19$ $e=3$<br>$G=A_4,\chi=\epsilon_7$ | $N=7$ $\epsilon=\epsilon_7$<br>$x^4+3x^2-7x+4$ | $F(10,5,0)$, $F(16,11,6)$, $F(22,17,12)$<br>$F(22,11,0)$, $F(28,17,6)$, $F(34,23,12)$ |
| $p=19$ $e=3$<br>$G=A_4,\chi=\epsilon_{11}$ | $N=11$ $\epsilon=\epsilon_{11}$<br>$x^4+15x^2-11x+81$ | $F(10,5,0)$, $F(16,11,6)$, $F(22,17,12)$<br>$F(22,11,0)$, $F(28,17,6)$, $F(34,23,12)$ |
| $p=7$ $e=3$<br>$G=S_4,\chi=1$ | $N=53$ $\epsilon=\epsilon_{53}$<br>$x^4-x^3+4x^2+1$ | $F(2,1,0)$, $F(4,3,2)$, $F(6,5,4)$<br>$F(6,3,0)$, $F(8,5,2)$, $F(10,7,4)$ |
| $p=13$ $e=4$<br>$G=S_4,\chi=1$ | $N=19$ $\epsilon=\epsilon_{19}$<br>$x^4-x^3+2x^2+4x-88$ | $F(7,5,3)$, $F(13,8,6)$, $F(16,14,9)$<br>$F(16,8,3)$, $F(19,14,6)$, $F(25,17,9)$ |
| $p=7$ $e=3$<br>$G=A_5,\chi=\epsilon_{73}$ | $N=73$ $\epsilon=\epsilon_{73}$<br>$x^5-5x^3-x^2+9x+7$ | $F(2,1,0)$, $F(4,3,2)$, $F(6,5,4)$<br>$F(6,3,0)$, $F(8,5,2)$, $F(10,7,4)$ |

TABLE 4. Irreducible higher level niveau 1 examples

5.3.2. *Representations of type $S_4$.* For $p>3$, $S_4$ has two absolutely irreducible three-dimensional representations defined over $\mathbb{F}_p$. Hence, by finding extensions $K/\mathbb{Q}$ with Galois group $S_4$ we may easily construct irreducible three-dimensional Galois representations which have image isomorphic to $S_4$. Two such examples are given in Table 4. Here, the format is as in the $A_4$ case, except that we take $\rho$ to be the unique irreducible three-dimensional representation of $G_{\mathbb{Q}}$ cutting out the splitting field of $f$, and taking transpositions to elements of trace 1 (twisted by $\chi$, which in both cases is trivial).

5.3.3. *Representations of type $A_5$.* The group $A_5$ has two three-dimensional irreducible representations defined over $\bar{\mathbb{F}}_p$, for each $p>5$. By composing these representations with the projection $G_{\mathbb{Q}}\to\mathrm{Gal}(K/\mathbb{Q})$, where $K$ is a field with Galois group $A_5$, we obtain irreducible three-dimensional Galois representations with image isomorphic to $A_5$. We give one example in Table 4, which we explain in detail.

*Example* 5.5. Let $K$ be the splitting field of the polynomial $x^5-5x^3-x^2+9x+7$. Then $\mathrm{Gal}(K/\mathbb{Q})$ is isomorphic to $A_5$, and $K$ is ramified only at $p=7$ (with ramification index 3) and at 73 (with ramification index 2). Let $\rho_1$ and $\rho_1'$ be the two characteristic 7 Galois representations alluded to above. Then it is easy to see that $\rho_1$ and $\rho_1'$ are Galois conjugates of each other over the field $\mathbb{F}_7$. The trace of both $\rho_1$ and $\rho_1'$ on a generator of inertia at 73 is $-1$, so that both representations have level $73^2$ and trivial nebentype. This level is too large for us to work with, so we twist both representations by the character $\chi=\epsilon_{73}$, to obtain $\rho=\rho_1\otimes\epsilon_{73}$, and $\rho'=\rho_1'\otimes\epsilon_{73}$. Now $\rho$ and $\rho'$ have level 73 and nebentype $\epsilon_{73}$.

Just as in Example 5.4, the restriction of $\rho$ (and of $\rho'$) to inertia at 7 has diagonal characters $\omega^0$, $\omega^2$, and $\omega^4$. Hence the predicted weights are the same as in those examples, namely $F(2,1,0)\otimes\det^a$ and $F(6,3,0)\otimes\det^a$ with $a=0,2,4$.

Computing the cohomology in each of these six weights with level 73 and nebentype $\epsilon_{73}$, we find that there is a unique eigenspace with the correct eigenvalues to correspond to $\rho$, and a unique eigenspace with the correct eigenvalues to correspond to $\rho'$ (at least up to $\ell=47$). As expected, these eigenspaces are defined over $\mathbb{F}_{7^2}$, rather than over $\mathbb{F}_7$, and they are Galois conjugates of each other over $\mathbb{F}_7$.

5.3.4. *Wildly ramified representations.* In addition to the preceding representations, we are able to calculate cohomology classes corresponding to irreducible three-dimensional representations $\rho : G_{\mathbb{Q}} \to \mathrm{GL}_3(\bar{\mathbb{F}}_p)$ which are wildly ramified at $p$. We have two types of examples of such representations; those having image $A_5$ which are wildly ramified at 5, and those having image $\mathrm{PSL}_2(\mathbb{F}_7)$, which are wildly ramified at 7.

We begin our study of the type $A_5$ representations by noting that there is a unique (up to isomorphism) injective homomorphism from $A_5$ to $\mathrm{GL}_3(\mathbb{F}_5)$, with image generated by the three matrices

$$\begin{pmatrix}1&1&0\\&1&1\\&&1\end{pmatrix},\quad\begin{pmatrix}4&2&2\\&1&2\\&&4\end{pmatrix},\quad\begin{pmatrix}4&1&4\\0&4&1\\2&4&2\end{pmatrix},$$

of orders 5, 2, and 3, respectively. The fields from which we obtain our Galois representations will have inertia group at 5 of order 5 or 10.

In the case of representations with inertia group of order 10, we will choose our representation so that the image of inertia is generated by the first two matrices above. With this choice of Galois representation it is clear that we have

$$\rho|_{I_p}\sim\begin{pmatrix}\omega^2&*&*\\&\omega^0&*\\&&\omega^2\end{pmatrix}.$$

Hence, we obtain a triple of $(2,0,2)$ yielding a predicted weight of

$$F(0,-1,2)'=F(4,3,2)=F(2,1,0)\otimes\det^2.$$

In order to keep the level to a manageable size, we will work with a twist of $\rho$ by a character unramified at $p$ (so that the weight is not affected). Let $\epsilon$ be the product of the characters $\epsilon_q$, where $q$ runs through the set of primes at which $\rho$ is ramified with ramification index 2. Then each prime $q$ at which $\rho$ has ramification index 2 contributes a factor of $q$ to the level of $\rho\otimes\epsilon$, and each prime $q$ at which $\rho$ has ramification index 3 contributes a factor of $q^2$ to the level of $\rho\otimes\epsilon$. The nebentype of $\rho\otimes\epsilon$ is easily seen to be $\epsilon$.

We have one example in which the inertia group has order 5. In this case we choose the representation so that the image of inertia is generated by the first matrix above. It is then clear that

$$\rho|_{I_p}\sim\begin{pmatrix}1&*&*\\&1&*\\&&1\end{pmatrix},$$

yielding a predicted weight of $F(-2,-1,0)'=F(6,3,0)$. The level of this representation is $3^4$ (note that 3 is wildly ramified) and the nebentype is trivial. For all of these examples, we have found that the predicted eigenclass does exist in the given weight, level and character, and has the correct eigenvalues (at least up to $\ell = 47$) to have $\rho$ attached.

| Polynomial | $G$ | $p$ | Weight | Level | $\epsilon$ |
|---|---|---|---|---|---|
| $x^5+5x^3-10x^2-45$ | $A_5$ | 5 | $F(4,3,2)$ | 13 | $\epsilon_{13}$ |
| $x^5+5x^3-10x^2-1$ | $A_5$ | 5 | $F(4,3,2)$ | 31 | $\epsilon_{31}$ |
| $x^5+5x^3-10x^2+5$ | $A_5$ | 5 | $F(4,3,2)$ | 37 | $\epsilon_{37}$ |
| $x^5+5x^3-10x^2+9$ | $A_5$ | 5 | $F(4,3,2)$ | 41 | $\epsilon_{41}$ |
| $x^5+5x^3-10x^2+20$ | $A_5$ | 5 | $F(4,3,2)$ | $2^2\cdot 13$ | $\epsilon_{13}$ |
| $x^5+25x^2+75$ | $A_5$ | 5 | $F(6,3,0)$ | $3^4$ | 1 |
| $x^7-7x^5-7x^4-7x^3-7x^2-7$ | $PSL_2(\mathbb{F}_7)$ | 7 | $F(6,5,4)$ | 17 | $\epsilon_{17}$ |
| $x^7+14x^6+14x^5-14x^4+35$ | $PSL_2(\mathbb{F}_7)$ | 7 | $F(8,5,2)$ | $5^2$ | 1 |
| $x^7-21x^3+7x-27$ | $PSL_2(\mathbb{F}_7)$ | 7 | $F(6,5,4)$ | 47 | $\epsilon_{47}$ |
| $x^7-7x^5-28x^2+7x+4$ | $PSL_2(\mathbb{F}_7)$ | 7 | $F(8,5,2)$ | $2^6$ | 1 |
| $x^7-7x^5-21x^4-49x^3-21x^2+1$ | $PSL_2(\mathbb{F}_7)$ | 7 | $F(8,5,2)$ | $2^6$ | 1 |
| $x^7-14x^4+42x^2-21x-9$ | $PSL_2(\mathbb{F}_7)$ | 7 | $F(6,5,4)$ | $3^4$ | 1 |
| $x^7+7x^5-7x^4-49x^3-98x-107$ | $PSL_2(\mathbb{F}_7)$ | 7 | $F(6,5,4)$ | $11^2$ | 1 |

TABLE 5. Wildly ramified Galois representations in niveau 1

We have also found Galois representations with image isomorphic to $\mathrm{PSL}_2(\mathbb{F}_7)$ which have low enough level that we can compute the predicted cohomology classes. The image of the representation is generated by the three matrices

$$\begin{pmatrix}1&1&0\\0&1&1\\0&0&1\end{pmatrix}, \quad \begin{pmatrix}2&3&3\\0&1&3\\0&0&4\end{pmatrix}, \quad \begin{pmatrix}1&0&0\\2&6&0\\4&0&6\end{pmatrix},$$

of orders 7, 3, and 4, respectively. The image of inertia at 7 under the representations that we will examine always has order 21, and may be chosen to be the subgroup generated by the first two matrices above. Hence, on inertia, we have

$$\rho|_{I_7} \sim \begin{pmatrix}\omega^2&*&*\\&\omega^0&*\\&&\omega^4\end{pmatrix}, \quad \text{or} \quad \begin{pmatrix}\omega^4&*&*\\&\omega^0&*\\&&\omega^2\end{pmatrix}.$$

In order to distinguish between these two possibilities, we use the action of tame ramification on wild ramification. Let $K$ be our $\mathrm{PSL}_2(\mathbb{F}_7)$ extension, and let $K_\mathfrak{p}$ be its localization at a prime above 7. Then $K_\mathfrak{p}/\mathbb{Q}_p$ is a degree 21 extension, which is totally ramified. Denote its Galois (inertia) group by $G_0$, and its higher ramification subgroups by $G_1, G_2, \ldots$. Clearly, there is a unique $i \geq 1$ such that $G_i/G_{i+1}$ is nontrivial, since $G_1$ is simple. By [20], for $a \in G_0/G_1$ and $b \in G_j/G_{j+1}$, we have

$$\theta_j(aba^{-1}) = \theta_0(a)^j\theta_j(t),$$

where $\theta_j : G_j/G_{j+1} \to U_K^{(j)}/U_K^{(j+1)}$ is the injective homomorphism sending $\sigma$ to $\sigma(\pi)/\pi$, where $\pi$ is a uniformizer of $K_\mathfrak{p}$.

If we consider $\theta_0$ as a map into $\mathbb{F}_p$, then we see (by [18]) that $\theta_0 = \omega^2$. We identify $G_i/G_{i+1} \cong G_i$ with its image under $\theta_i$. Then we have

$$sts^{-1} = t^{\omega^{2i}(s)}.$$

| Fixed field of $\ker(\sigma)$ | $G$ | $k$ | $p$ | Weights | | Level | $\epsilon$ |
|---|---|---|---|---|---|---|---|
| $x^3 - x^2 - 8x + 7$ | $S_3$ | 1 | 5 | $F(5,4,1)$ | $F(4,4,2)$ | 73 | $\epsilon_{73}$ |
| | | 3 | 5 | $F(4,2,2)$ | $F(5,2,1)$ | 73 | $\epsilon_{73}$ |
| $x^3 - x^2 - 7x + 2$ | $S_3$ | 5 | 11 | $F(5,4,3)$ | $F(22,14,6)$ | 13 | $\epsilon_{13}$ |
| | | 7 | 11 | $F(15,6,3)$ | $F(12,6,6)$ | 13 | $\epsilon_{13}$ |
| | | 9 | 11 | $F(15,8,3)$ | $F(12,8,6)$ | 13 | $\epsilon_{13}$ |
| | | 11 | 11 | $F(15,10,3)$ | $F(12,10,6)$ | 13 | $\epsilon_{13}$ |
| | | 13 | 11 | $F(15,12,3)$ | $F(12,12,6)$ | 13 | $\epsilon_{13}$ |
| $x^3 - 11x - 11$ | $S_3$ | 5 | 11 | $F(5,4,3)$ | $F(22,14,6)$ | 17 | $\epsilon_{17}$ |
| | | 7 | 11 | $F(15,6,3)$ | $F(12,6,6)$ | 17 | $\epsilon_{17}$ |
| | | 9 | 11 | $F(15,8,3)$ | $F(12,8,6)$ | 17 | $\epsilon_{17}$ |
| | | 11 | 11 | $F(15,10,3)$ | $F(12,10,6)$ | 17 | $\epsilon_{17}$ |
| | | 13 | 11 | $F(15,12,3)$ | $F(12,12,6)$ | 17 | $\epsilon_{17}$ |

TABLE 6. Reducible niveau 2 representations $\sigma \oplus \omega^k$ with $\sigma$ even

However, we see easily (by matrix multiplication) that $sts^{-1} = t^2$, so that $\omega^{2i}(s) = 2$, and we have that

$$\rho|_{I_7} \sim \begin{pmatrix} \omega^{2i} & * & * \\ & 1 & * \\ & & \omega^{4i} \end{pmatrix}.$$

Finally, analysis of the ramification groups shows that if the discriminant of the degree 7 subfield of $K$ is exactly divisible by $7^8$ then $i = 1$, and if it is exactly divisible by $7^{10}$ then $i = 2$.

Clearly if $i = 1$, we get a predicted weight of $F(0,-1,4)' = F(6,5,4)$, and if $i = 2$, we get a predicted weight of $F(2,-1,2)' = F(8,5,2)$. The level and nebentype are easily calculated, and in case of odd primes $q$ which have inertia group of order 2, we twist by $\epsilon_q$ to lower the level from $q^2$ to $q$.

Table 5 contains information on the wildly ramified Galois representations which we have studied. Each line of the table gives a polynomial whose Galois closure is a $G$-extension of $\mathbb{Q}$ (where $G = A_5$ or $\mathrm{PSL}_2(\mathbb{F}_7)$), as well as the weight, level and nebentype predicted by Conjecture 3.1 for the appropriate twist of $\rho$. In each case, the representation for which the invariants were computed is actually $\rho \otimes \epsilon$, where $\epsilon$ is the indicated nebentype (as described above, this lowers the level to a manageable size). In every example an eigenclass with the correct eigenvalues (up to $\ell = 47$) occurs in the predicted cohomology group.

## 6. Niveau 2 representations

6.1. **Reducible representations in higher level.** Each line of Table 6 gives a polynomial with splitting field a totally real $S_3$-extension $K$ of $\mathbb{Q}$. In each case, we define $\sigma$ to be the unique two-dimensional Galois representation $\sigma : G_{\mathbb{Q}} \to \mathrm{GL}_2(\bar{\mathbb{F}}_p)$ which cuts out $K$, and we note that $\sigma$ is niveau 2. Letting $\rho = \sigma \oplus \omega^k$, where $k$ is indicated in the table, Conjecture 3.1 predicts two possible weights corresponding to $\rho$, as indicated in the table. We have checked that for each row of Table 6 the cohomology in the given weights, level, and nebentype does contain an eigenclass with the correct eigenvalues to correspond to $\rho$, at least up to $\ell = 47$.

Table 7 contains examples of Galois representations each of which is the sum of an odd irreducible two-dimensional Galois representation and the trivial character.

| Galois Representation | | Predicted Weights |
|---|---|---|
| $p = 5$, $e = 3$<br>$G = S_3$ | $N = 7$ $\epsilon = \epsilon_3$<br>$x^3 - x^2 + 2x - 3$ | $F(1,0,0)$, $F(4,1,0)$<br>$F(2,2,1)$, $F(6,5,2)$ |
| $p = 5$, $e = 3$<br>$G = S_3$ | $N = 43$ $\epsilon = \epsilon_3$<br>$x^3 - x^2 + 2x + 12$ | $F(1,0,0)$, $F(4,1,0)$<br>$F(2,2,1)$, $F(6,5,2)$ |
| $p = 5$, $e = 3$<br>$G = S_3$ | $N = 47$ $\epsilon = \epsilon_3$<br>$x^3 + 5x - 5$ | $F(1,0,0)$, $F(4,1,0)$<br>$F(2,2,1)$, $F(6,5,2)$ |
| $p = 5$, $e = 3$<br>$G = S_3$ | $N = 67$ $\epsilon = \epsilon_3$<br>$x^3 - x^2 + 7x + 2$ | $F(1,0,0)$, $F(4,1,0)$<br>$F(2,2,1)$, $F(6,5,2)$ |
| $p = 5$, $e = 3$<br>$G = S_3$ | $N = 83$ $\epsilon = \epsilon_3$<br>$x^3 - 10x - 15$ | $F(1,0,0)$, $F(4,1,0)$<br>$F(2,2,1)$, $F(6,5,2)$ |
| $p = 5$, $e = 3$<br>$G = S_3$ | $N = 83$ $\epsilon = \epsilon_3$<br>$x^3 - x^2 + 7x - 8$ | $F(1,0,0)$, $F(4,1,0)$<br>$F(2,2,1)$, $F(6,5,2)$ |
| $p = 5$, $e = 3$<br>$G = S_3$ | $N = 83$ $\epsilon = \epsilon_3$<br>$x^3 - x^2 - 3x - 8$ | $F(1,0,0)$, $F(4,1,0)$<br>$F(2,2,1)$, $F(6,5,2)$ |
| $p = 17$, $e = 3$<br>$G = S_3$ | $N = 3$ $\epsilon = \epsilon_3$<br>$x^3 - x^2 + 6x - 12$ | $F(9,4,0)$, $F(20,9,0)$<br>$F(14,10,5)$, $F(30,21,10)$ |
| $p = 17$, $e = 3$<br>$G = S_3$ | $N = 7$ $\epsilon = \epsilon_7$<br>$x^3 - x^2 + 6x + 5$ | $F(9,4,0)$, $F(20,9,0)$<br>$F(14,10,5)$, $F(30,21,10)$ |
| $p = 3$, $e = 4$<br>$G = D_4$ | $N = 7$ $\epsilon = \epsilon_7$<br>$x^4 - 3x^2 - 3$ | $F(1,0,0)$, $F(2,2,1)$<br>$F(2,1,0)$ |
| $p = 3$, $e = 4$<br>$G = D_4$ | $N = 19$ $\epsilon = \epsilon_{19}$<br>$x^4 - 30x^2 - 3$ | $F(1,0,0)$, $F(2,2,1)$<br>$F(2,1,0)$ |
| $p = 3$, $e = 4$<br>$G = D_4$ | $N = 31$ $\epsilon = \epsilon_{31}$<br>$x^4 + 9x^2 - 3$ | $F(1,0,0)$, $F(2,2,1)$<br>$F(2,1,0)$ |
| $p = 3$, $e = 4$<br>$G = D_4$ | $N = 43$ $\epsilon = \epsilon_{43}$<br>$x^4 - 318x^2 - 3$ | $F(1,0,0)$, $F(2,2,1)$<br>$F(2,1,0)$ |
| $p = 7$, $e = 4$<br>$G = D_4$ | $N = 11$ $\epsilon = \epsilon_{11}$<br>$x^4 - 7x^2 - 7$ | $F(3,0,0)$, $F(6,3,0)$<br>$F(4,4,1)$, $F(10,7,4)$ |
| $p = 19$, $e = 5$<br>$G = D_5$ | $N = 7$, $\epsilon = \epsilon_7$<br>$x^5 - 19x^2 + 38x - 95$ | $F(13,2,0)$, $F(20,13,0)$<br>$F(16,14,3)$, $F(34,21,14)$ |
| $p = 19$, $e = 5$<br>$G = D_5$ | $N = 7$, $\epsilon = \epsilon_7$<br>$x^5 - 19x^2 + 38x - 95$ | $F(9,6,0)$, $F(24,9,0)$<br>$F(16,10,7)$, $F(34,25,10)$ |

TABLE 7. Reducible niveau 2 representations $\sigma \oplus \omega^0$ with $\sigma$ odd

In each example, a polynomial is given which has Galois group $G$. For all but the last two examples we let $\sigma$ be the unique two-dimensional mod $p$ representation of $G$, and in all cases we take $\rho$ to be $\sigma \oplus 1$. The ramification index $e$ of $p$, and the level $N$ and nebentype $\epsilon$ of $\rho$ are indicated. For each such representation Conjecture 3.1 predicts four weights (two of the predicted weights are the same in the $p = 3$ cases), as indicated in the table, and in all cases we have been able to check that the predicted eigenvalues exist in the cohomology in all of the predicted weights. We explain the last two examples in Table 7 in detail in Example 6.1.

*Example* 6.1. Let $K$ be the splitting field of the polynomial $f = x^5 - 19x^2 + 38x - 95$. Then $K$ is a totally complex $D_5$-extension of $\mathbb{Q}$, ramified only at 7 (with ramification index 2 and residual degree 1) and 19 (with ramification index 5 and residual degree 2).

The group $D_5$ has two irreducible two-dimensional representations—we will denote the corresponding Galois representations by $\sigma$ and $\sigma'$. Let $\rho$ (resp. $\rho'$) be the direct sum of $\sigma$ (resp. $\sigma'$) with the trivial character. Then both $\rho$ and $\rho'$ are easily seen to have level 7 and nebentype $\epsilon_7$.

We may conjugate each of $\rho$ and $\rho'$ to land in either of the standard Levi subgroups

$$L = \begin{pmatrix} * & * & \\ * & * & \\ & & * \end{pmatrix}, \quad or \quad L' = \begin{pmatrix} * & & \\ & * & * \\ & * & * \end{pmatrix},$$

and each representation will satisfy strict parity with Levi subgroup $L$ (or $L'$), as $\sigma$ and $\sigma'$ are both odd.

Both $\rho$ and $\rho'$ have niveau two, but they differ on restriction to inertia at 19. We choose $\rho$ (possibly swapping $\sigma$ and $\sigma'$) so that

$$\rho|_{I_{19}} \sim_L \begin{pmatrix} \psi^{72} & & \\ & {\psi'}^{72} & \\ & & \omega^0 \end{pmatrix},$$

while $\rho'$ will have diagonal characters $\psi^{144}, {\psi'}^{144}, \omega^0$.

Since $72 = 15+3*19$, we get a predicted weight of $F(15-2, 3-1, 0)' = F(13, 2, 0)$ for $\rho$. In addition, we may also conjugate $\rho$ inside $L$ so that the diagonal characters on inertia are $\psi^{288}, {\psi'}^{288}$ and $\omega^0$. Since $288 = 22 + 14 * 19$ we also predict a weight of $F(20, 13, 0)$ for $\rho$. Finally we may conjugate $\rho$ to land inside $L'$ yielding predicted weights of $F(16, 14, 3)$ and $F(34, 21, 14)$. In a similar fashion, we predict four weights for $\rho'$, namely $F(9, 6, 0), F(16, 10, 7), F(24, 9, 0)$, and $F(34, 25, 10)$.

In order to test whether the representations $\rho$ and $\rho'$ are attached to Hecke eigenclasses with these weights, we need to compute the characteristic polynomials of the images of Frobenius elements under $\rho$ and $\rho'$. There is a subtlety introduced here by the fact that $D_5$ has two conjugacy classes of order 5. On one of these classes $\rho$ has trace 5 and $\rho'$ has trace 15, while on the other class these traces are reversed. We must determine which class contains each Frobenius element of order 5.

Suppose $\tau \in G_{\mathbb{Q}}$ restricts to an element of order 5 in $\mathrm{Gal}(K/\mathbb{Q})$. Then there is some element $\eta \in I_5$ such that $\tau \equiv \eta$ modulo $G_K$. So

$$\mathrm{Tr}(\rho(\tau)) = \mathrm{Tr}(\rho(\eta)) = \psi^{72}(\eta) + {\psi'}^{72}(\eta) + 1.$$

Let $\mathfrak{P}$ be the unique prime of $K$ lying over $p = 19$, and let $\pi$ be a uniformizer of $\mathfrak{P}$. Then $\zeta = \psi^{72}(\eta) \equiv \eta(\pi)/\pi \pmod{\mathfrak{P}}$ is a fifth root of unity in the residue field $F$ of $\mathfrak{P}$ (which has order $19^2$). Note that there are two possible images of $\zeta$ in $F$, depending on our choice of fundamental character $\psi$. We may, however, compute $\zeta + \zeta^p$, which will be in the prime field, and is independent of this choice. These calculations are easily performed using GP/PARI, since $K$ is only of degree 10 over $\mathbb{Q}$. A convenient uniformizer to use is a root $\alpha$ of the polynomial $f$ defined above.

We find, for instance, that $\mathrm{Tr}(\rho(\mathrm{Frob}_2)) = 5$, giving predicted Hecke eigenvalues $a(2, 1) = 5$, and $a(2, 2) = 12$ for the classes attached to $\rho$, and eigenvalues $a(2, 1) = 15$ and $a(2, 2) = 17$ for the classes attached to $\rho'$.

We have computed the Hecke eigenvalues (for $l \leq 47$) for cohomology classes with each of the 8 weights that arose above, and in each case have found that the eigenvalues are exactly as predicted.

6.2. **Irreducible representations in higher level.** In niveau 2, we again obtain several irreducible representations which are symmetric squares of odd two-dimensional representations, but we also obtain one set of examples which are not. We begin by describing an example of the former type of representation.

*Example* 6.2. Let $p = 5$ and let $K$ be the splitting field of the polynomial $f = x^4 + x^2 - x + 2$. Then $K$ has Galois group $S_4$, and is ramified only at 5 (with ramification index 3) and at 73 (with ramification index 2). In fact, since the discriminant of $f$ is $5^2 73$, the quadratic subfield of $K$ is ramified at 73, so the inertia group at 73 must be generated by a transposition. If we let $\rho$ be the unique irreducible three-dimensional mod 5 representation of $G_{\mathbb{Q}}$ cutting out $K$ and taking transpositions to elements of trace 1, we easily determine that the level of $\rho$ is 73, and that its nebentype is $\epsilon_{73}$. The weights predicted for $\rho$ by Conjecture 3.1 are calculated by noting that

$$\rho|_{I_5} \sim \begin{pmatrix} \psi^8 & & \\ & {\psi'}^8 & \\ & & \omega^0 \end{pmatrix},$$

with $8 = 3{+}1{*}5$, so that we have a predicted weight of $F(3{-}2, 1{-}1, 0)' = F(1, 0, 0)'$. We may also write

$$\rho|_{I_5} \sim \begin{pmatrix} \psi^8 & & \\ & \omega^0 & \\ & & {\psi'}^8 \end{pmatrix}, \quad \text{or} \quad \rho|_{I_5} \sim \begin{pmatrix} \omega^0 & & \\ & \psi^8 & \\ & & {\psi'}^8 \end{pmatrix},$$

yielding predicted weights of $F(3-2, 0-1, 1)' = F(5, 3, 1)$ and $F(0-2, 3-1, 1)' = F(2, 2, 1)'$.

In addition, we note that ${\psi'}^8 = \psi^{16}$, so we can permute the two characters of niveau 2, and write

$$\rho|_{I_5} \sim \begin{pmatrix} \psi^{16} & & \\ & {\psi'}^{16} & \\ & & \omega^0 \end{pmatrix},$$

with $16 = 1{+}3{*}5 = 6{+}2{*}5$, so that we have a predicted weight of $F(6{-}2, 2{-}1, 0)' = F(4, 1, 0)$. We may also write

$$\rho|_{I_5} \sim \begin{pmatrix} \psi^{16} & & \\ & \omega^0 & \\ & & {\psi'}^{16} \end{pmatrix}, \quad \text{or} \quad \rho|_{I_5} \sim \begin{pmatrix} \omega^0 & & \\ & \psi^{16} & \\ & & {\psi'}^{16} \end{pmatrix},$$

yielding predicted weights of $F(6-2, 0-1, 2)' = F(4, 3, 2)$ and $F(0-2, 6-1, 2)' = F(6, 5, 2)$.

Calculating the cohomology in all six of these weights, we find eigenclasses with the correct Hecke eigenvalues to correspond to $\rho$ (at least for primes up to 47). This yields a family of six "companion forms" of different weights, all of which seem to correspond to $\rho$.

Other examples of irreducible niveau 2 representations with image isomorphic to $S_4$, as well as examples with image isomorphic to $A_5$ (where the representation is the twist by $\chi$ of the unique irreducible three-dimensional mod 5 representation having image $A_5$ and cutting out the splitting field of $f$), are given in Table 8, in the same format as the examples in Table 4. In addition, examples with image of order 54 are given in [12]. These last examples have $p = 5$, level $N = 83$, with quadratic nebentype, and can not be the symmetric square of any two-dimensional representation.

| Galois Representation | | Predicted Weights |
|---|---|---|
| $p = 5$ $e = 3$<br>$G = S_4, \chi = 1$ | $N = 73$ $\epsilon = \epsilon_{73}$<br>$x^4 + x^2 - x + 2$ | $F(1,0,0)$, $F(5,3,1)$, $F(2,2,1)$<br>$F(4,1,0)$, $F(4,3,2)$, $F(6,5,2)$ |
| $p = 5$ $e = 3$<br>$G = S_4$, $\chi = 1$ | $N = 144$ $\epsilon = 1$<br>$x^4 - 2x^3 - 8x + 4$ | $F(1,0,0)$, $F(5,3,1)$, $F(2,2,1)$<br>$F(4,1,0)$, $F(4,3,2)$, $F(6,5,2)$ |
| $p = 7$, $e = 4$<br>$G = S_4$, $\chi = 1$ | $N = 67$ $\epsilon = \epsilon_{67}$<br>$x^4 - 56x + 112$ | $F(6,3,3)$, $F(12,8,4)$, $F(7,7,4)$<br>$F(9,6,3)$, $F(3,2,1)$, $F(7,4,1)$ |
| $p = 11$, $e = 3$<br>$G = S_4$, $\chi = 1$ | $N = 17$ $\epsilon = \epsilon_{17}$<br>$x^4 - x^3 + 3x + 2$ | $F(5,2,0)$, $F(15,9,3)$, $F(8,6,3)$<br>$F(12,5,0)$, $F(12,9,6)$, $F(18,13,6)$ |
| $p = 5$, $e = 3$<br>$G = A_5$, $\chi = \epsilon_{89}$ | $N = 89$ $\epsilon = \epsilon_{89}$<br>$x^5 - 2x^3 - x^2 - 6x - 11$ | $F(1,0,0)$, $F(5,3,1)$, $F(2,2,1)$<br>$F(4,1,0)$, $F(4,3,2)$, $F(6,5,2)$ |
| $p = 5$, $e = 3$<br>$G = A_5$, $\chi = \epsilon_{151}$ | $N = 151$, $\epsilon = \epsilon_{151}$<br>$x^5 - 3x^3 - x^2 + x - 3$ | $F(1,0,0)$, $F(5,3,1)$, $F(2,2,1)$<br>$F(4,1,0)$, $F(4,3,2)$, $F(6,5,2)$ |
| $p = 5$, $e = 3$<br>$G = A_5$, $\chi = \epsilon_{157}$ | $N = 157$ $\epsilon = \epsilon_{157}$<br>$x^5 + 7x^3 - x^2 - 9x + 7$ | $F(1,0,0)$, $F(5,3,1)$, $F(2,2,1)$<br>$F(4,1,0)$, $F(4,3,2)$, $F(6,5,2)$ |

TABLE 8. Irreducible niveau 2 representations

## 7. Niveau 3 representations

We have two examples of odd niveau 3 representations, both of which support Conjecture 3.1. It is easy to see that a niveau 3 representation must be irreducible, and that it cannot be the symmetric square of a two-dimensional representation. Our first example is induced from a one-dimensional representation on a subgroup of index 3 in $G_{\mathbb{Q}}$, and the second has image isomorphic to $\mathrm{PSL}_2(\mathbb{F}_7)$ in $\mathrm{GL}_3(\bar{\mathbb{F}}_{11})$, but is in no obvious way related to any representation of dimension less than 3.

7.1. **An induced representation.** Let $f = x^3 + 2x - 1$. The Galois group of $f$ is $S_3$. Let $K$ be the splitting field of $f$, and $K_3 = \mathbb{Q}(\alpha)$, where $\alpha$ is a root of $f$. Then $K_3$ is ramified only at 59. Using GP/PARI, we may calculate the ray class group of $K_3$ modulo 7, and find that it is cyclic of order 9. If we let $L$ be the ray class field of $K_3$ modulo 7, then the existence of $L$ implies the existence of a character $\chi : G_{K_3} \to \mu_9 \subset \mathbb{F}_{7^3}$ of order 9, ramified only at primes above 7. If we now set

$$\rho = \mathrm{Ind}_{G_{K_3}}^{G_{\mathbb{Q}}} \chi,$$

then $\rho : G_{\mathbb{Q}} \to \mathrm{GL}_3(\bar{\mathbb{F}}_7)$ must be irreducible, since it has niveau 3 (as the ramification index at 7 is divisible by 9). Note that there are six choices of $\chi$, since there are six primitive ninth roots of unity in $\mu_9$. Until we make a choice, everything we state will be true for any choice of $\chi$, and hence for any $\rho$ induced from $\chi$.

If we let $M$ be the Galois closure of $L$, we see that $M$ contains the composite field $KL$, which is abelian of degree 9 over $K$, and in fact, $M$ is generated by the conjugates of $KL$ over $\mathbb{Q}$. We see from this that no element of $\mathrm{Gal}(M/K)$ has order more than 9. Note that $\rho$ factors through $G = \mathrm{Gal}(M/\mathbb{Q})$, so in particular, the image of inertia at 7 under $\rho$ must be of order 9 (since inertia fixes $K$). In fact, it is easy to see that the factorization of $\rho$ through $G$ is a faithful representation of $G$.

Now let

$$G_{\mathbb{Q}} = \bigcup_{i=0}^{2} g_i G_{K_3},$$

where the $g_i$ are coset representatives of $G_{K_3}$ in $G_{\mathbb{Q}}$, and we have that for $g \in G_{\mathbb{Q}}$,

$$\mathrm{Tr}(\rho(g)) = \sum_{i=0}^{2} \chi^0(g^{g_i})$$

where

$$\chi^0(x) = \begin{cases} 0 & \text{if } x \notin G_{K_3} \\ \chi(x) & \text{if } x \in G_{K_3}. \end{cases}$$

Using this description of $\rho$, we may calculate values of $\mathrm{Tr}(\rho(g))$ in terms of $\chi$ for various $g$, given that we know the order of $\pi(g)$, where $\pi : G_{\mathbb{Q}} \to S_3$ is the natural projection onto the Galois group of $K$. Let $g' \in G_{K_3}$ be a conjugate by some $g_i$ of $g$, if such a conjugate exists, and we have that for $\pi(g)$ of order 2, $\rho(g)$ has trace $\chi(g')$, and for $\pi(g)$ of order 3, $\rho(g)$ has trace 0 (since no conjugate of $g$ is in $G_{K_3}$).

In fact, we may even go further and compute the values of $\chi(g')$ using class field theory. Class field theory shows the existence of an isomorphism between the ray class group $J$ of $K_3$ modulo 7 and the group $\mu_9$ of ninth roots of unity. We fix this isomorphism by setting the image of the ideal $\mathfrak{p}$ above 2 in $K_3$ with inertial degree 1 to have image $\chi(\mathrm{Frob}_{\mathfrak{p}}) = \zeta_9$, where $\mathrm{Frob}_{\mathfrak{p}}$ is a Frobenius above $\mathfrak{p}$. (note that $\mathfrak{p}$ has order 9 in the ray class group). Given any ideal of $K_3$, we may then compute its image in $J$ in terms of the image of the ideal above 2, and hence find the image of any Frobenius element under $\chi$. The ray class computations are easily done using GP/PARI, since the degree of $K_3$ is only 3.

Using these techniques, we find the following values:

| $p$ | 2 | 3 | 5 | 7 | 11 | 13 | 17 | 19 | 23 | 29 | 31 | 37 | 41 | 43 | 47 |
|---|---|---|---|---|---|---|---|---|---|---|---|---|---|---|---|
| $O(\mathrm{Frob}_p)$ | 18 | 9 | 9 | * | 18 | 6 | 9 | 9 | 18 | 3 | 18 | 18 | 3 | 2 | 18 |
| $O(\pi(\mathrm{Frob}_p))$ | 2 | 3 | 3 | * | 2 | 2 | 1 | 3 | 2 | 3 | 2 | 2 | 3 | 2 | 2 |
| $\chi(\mathrm{Frob}'_p)$ | $\zeta_9$ | * | * | * | $\zeta_9^5$ | $\zeta_9^6$ | $\zeta_9$ | * | $\zeta_9$ | * | $\zeta_9^2$ | $\zeta_9^7$ | * | 1 | $\zeta_9$ |
| $\mathrm{Tr}(\rho(\mathrm{Frob}_p))$ | $\zeta_9$ | 0 | 0 | * | $\zeta_9^5$ | $\zeta_9^6$ | 0 | 0 | $\zeta_9$ | 0 | $\zeta_9^2$ | $\zeta_9^7$ | 0 | 1 | $\zeta_9$ |

The only value which has not yet been explained is that the Frobenius at 17 has image of trace 0. This occurs since 17 splits completely in $K$ (hence also in $K_3$). Hence, there are three distinct conjugates $\mathrm{Frob}_{17}^{g_i}$ of $\mathrm{Frob}_{17}$, all in $G_{K_3}$, and their images under $\chi$ are $\zeta_9$, $\zeta_9^4$, and $\zeta_9^7$, so that the trace of $\rho(\mathrm{Frob}_{17})$ is 0.

Direct computation in the ray class group shows that if $p \leq 47$ is a rational prime with $\pi(\mathrm{Frob}_p)$ having order 2, and $g$ is any Frobenius element for $p$, then $\chi(g^2) = \chi(g')^2$. Since this is true for any conjugate of $g^2$, we have that $\mathrm{Tr}(\rho(g^2)) = 3\chi(g')^2 = 3\mathrm{Tr}(\rho(g))^2$. Using this fact, a simple computation (using Magma) shows that the eigenvalues of $\rho(g)$ must be $\xi$, $\xi$ and $-\xi$, where $\xi = \mathrm{Tr}(\rho(g))$. Hence, the characteristic polynomial $\det(1 - \rho(g)X)$ is equal to

$$1 - \xi X - \xi^2 X^2 + \xi^3 X^3.$$

In particular, we will use the fact that $\det \rho(g) = -\xi^3 = -(\mathrm{Tr}\rho(g))^3$.

We now compute the level and character of $\rho$. The prime 59 has ramification index 2 in the fixed field of $\rho$, and if $g$ is a generator of inertia at 59, then $\pi(g)$ has order 2 (since 59 has ramification index 2 in $K$). In addition, $\chi(g')$ must be simultaneously a ninth root and a square root of 1, hence equal to 1. Then $\mathrm{Tr}(\rho(g)) = \chi(g') = 1$, so the three eigenvalues of $g$ are 1, 1, and $-1$, and the level of $\rho$ must be 59, with nebentype $\epsilon_{59}$.

Finally, we calculate the predicted weights for $\rho$. These weights will in fact depend on the choice of $\chi$. We recall that the fundamental characters of niveau 3 will be denoted by $\theta$, $\theta'$ and $\theta''$. Since 7 has ramification index 9 in $M$, we know that $\rho$ must have niveau 3. In fact, we have that either

$$\rho|_{I_7} \sim \begin{pmatrix} \theta^{38} & & \\ & {\theta'}^{38} & \\ & & {\theta''}^{38} \end{pmatrix} \quad \text{or} \quad \rho|_{I_7} \sim \begin{pmatrix} \theta^{76} & & \\ & {\theta'}^{76} & \\ & & {\theta''}^{76} \end{pmatrix}.$$

Note that

$$\det \rho = \omega^{38}\epsilon_{59} = \omega^2\epsilon_{59}$$

in the first case, while

$$\det \rho = \omega^{76}\epsilon_{59} = \omega^4\epsilon_{59}$$

in the second. Thus, in order to obtain the first case, we choose $\chi$ (and hence $\zeta_9 = \chi(\mathrm{Frob}_2)$) so that

$$-\zeta_9^3 = -\mathrm{Tr}(\rho(\mathrm{Frob}_2))^3 = \det(\rho(\mathrm{Frob}_2)) = \omega^2(\mathrm{Frob}_2)\epsilon_{59}(\mathrm{Frob}_2) = -4,$$

and in order to get the second we choose $\chi$ (and hence $\zeta_9$) so that $-\zeta_9^3 = -2$.

Note that each of the two possibilities comes from three choices of $\chi$. Hence, we should expect three eigenclasses in each predicted weight—one for each choice of $\chi$.

Considering the first case, $m = 38 = 3 + 5 * 7 + 0 * 7^2$, so we get a triple $(a, b, c) = (5, 3, 0)$, yielding predicted weight

$$F(5 - 2, 3 - 1, 0) = F(3, 2, 0).$$

We may also permute the characters on the diagonal, which will have the effect of multiplying $m$ by 7 or $7^2$ modulo $7^3 - 1$, yielding predicted triples and weights as follows:

For $7 * m = 266 = 7 + 9 * 7 + 4 * 7^2$ we get predicted weight

$$F(9 - 2, 7 - 1, 4) = F(3, 2, 0) \otimes \det{}^4.$$

For $49 * m \equiv 152 = 5 + 7 * 7 + 2 * 7^2$, we get predicted weight

$$F(7 - 2, 5 - 1, 2) = F(3, 2, 0) \otimes \det{}^2.$$

We may similarly calculate weights for the second possibility, and find the following predicted weights:

$$F(3, 1, 0) \otimes \det{}^1, \quad F(3, 1, 0) \otimes \det{}^3, \quad \text{and} \quad F(3, 1, 0) \otimes \det{}^5,$$

Computations show that cohomology classes with the correct eigenvalues (up to $\ell = 47$) exist in all of these weights. In each weight there is a triple of eigenclasses, defined over $\mathbb{F}_{7^3}$ and conjugate over $\mathbb{F}_7$, each corresponding to a choice of $\chi$ as above.

| Class | 1 | 2 | 3 | 4 | 5 | 6 |
|---|---|---|---|---|---|---|
| Size | 1 | 21 | 56 | 42 | 24 | 24 |
| Order | 1 | 2 | 3 | 4 | 7 | 7 |
| $\chi_1$ | 1 | 1 | 1 | 1 | 1 | 1 |
| $\chi_2$ | 3 | -1 | 0 | 1 | $\alpha$ | $\bar{\alpha}$ |
| $\chi_3$ | 3 | -1 | 0 | 1 | $\bar{\alpha}$ | $\alpha$ |
| $\chi_4$ | 6 | 2 | 0 | 0 | -1 | -1 |
| $\chi_5$ | 7 | -1 | 1 | -1 | 0 | 0 |
| $\chi_6$ | 8 | 0 | -1 | 0 | 1 | 1 |

TABLE 9. Character Table of $\mathrm{PSL}_2(\mathbb{F}_7)$

7.2. **A representation with image $\mathbf{PSL}_2(\mathbb{F}_7)$.** We begin by noting that the irreducible polynomial $f_1 = x^7 - 11x^5 - 22x^4 + 33x^2 + 33x + 11$ has Galois group $PSL_2(\mathbb{F}_7)$, as reported by GP/PARI, KASH, and Magma. If $L = \mathbb{Q}(\alpha)$, where $\alpha$ is a root of $f$, we find that the discriminant of $L$ is $11^6 31^2$. Since 11 is tamely ramified, we may conclude that the ramification index of 11 in the splitting field $K$ of $f$ is $e = 7$. Using the main result of [8], we see easily that the ramification index of 31 in $K$ is 2.

The character table of $PSL_2(\mathbb{F}_7)$ is given in Table 9, where $\alpha = \frac{-1+\sqrt{-7}}{2}$ and $\bar{\alpha} = \frac{-1-\sqrt{-7}}{2}$. Over $\bar{\mathbb{F}}_{11}$, we have that $\alpha$ and $\bar{\alpha}$ are equal to 4 and 6, with the order depending on our choice of $\sqrt{-7}$.

The existence of the $PSL_2(\mathbb{F}_7)$-extension $K$ gives rise to two irreducible three-dimensional Galois representations defined over $\bar{\mathbb{F}}_{11}$. The image of inertia at 11 under both representations has order 7, so they are both niveau 3. We choose $\sigma$ to be the representation which when restricted to inertia at 11 has diagonal characters $\theta^{190}$, $\theta'^{190}$, and $\theta''^{190}$, and $\sigma'$ to be the other (with diagonal characters on inertia equal to $\theta^{570}, \theta'^{570}, \theta''^{570}$).

We note that the level of $\sigma$ (and of $\sigma'$) is $31^2$, since the elements of order 2 are mapped to matrices of trace $-1$. This level is too large for convenient calculation, so we investigate $\rho = \sigma \otimes \epsilon_{31}$ and $\rho' = \sigma' \otimes \epsilon_{31}$, which are easily seen to have level 31 and nebentype $\epsilon_{31}$.

In order to calculate the predicted eigenvalues of the image of a Frobenius element of order 7 under $\rho$, we need to distinguish between the two conjugacy classes of order 7 in $\mathrm{PSL}_2(\mathbb{F}_7)$. In order to do this, we use a method similar to that used in Example 6.1. In this case, the method needs to be modified slightly, since we are dealing with much larger fields.

We begin by using Magma to determine the Galois group $G \cong \mathrm{PSL}_2(\mathbb{F}_7)$ of $f$ as a permutation group acting on the roots of $f$. We note that each root of $f$ is a uniformizer for all primes lying above 11 in $K$ (since 11 is tamely ramified, and all the ramification occurs in $L/\mathbb{Q}$). Let $\alpha$ be a root of $f$, and let $\tau$ be an element of order 7 in $G$. Then we easily compute a complex approximation to $\beta = \tau(\alpha)/\alpha$. If $\mathfrak{P}$ is the prime of $K$ lying over 11 and having inertia group generated by $\tau$, then the image of $\beta$ in the residue field of $\mathfrak{P}$ is a Galois conjugate of the primitive seventh root of unity $\theta^{190}(\tau)$. Hence, the trace of $\sigma(\tau)$ is equal to $\beta + \beta^{11} + \beta^{121}$ mod $\mathfrak{P}$. We actually compute a complex approximation to $\gamma = \beta + \beta^2 + \beta^4$, which is equal to this trace modulo $\mathfrak{P}$. Knowing that this trace is congruent to either 4 or 6 modulo $\mathfrak{P}$, we compute $\delta_1 = \gamma - 4$, and $\delta_2 = \gamma - 6$. Exactly one of $\delta_1$ and $\delta_2$ should lie in $\mathfrak{P}$. We note that if $K_8$ is the unique degree 8 subfield of $K$ fixed by $\langle\tau\rangle$ (so that $K_8$ is the decomposition field of $\mathfrak{P}$), then there is a unique degree 1 prime $\mathfrak{p}$ in $K_8$, and $\mathfrak{P}$ lies over $\mathfrak{p}$. Hence, we may determine whether $\delta_i$ lies in $\mathfrak{P}$ by determining whether the norm of $\delta_i$ (from $K$ to $K_8$) lies in $\mathfrak{p}$. We compute a complex approximation to this norm (and all of its Galois conjugates), and then easily find a complex approximation to the minimal polynomial of this norm. This polynomial should have rational integer coefficients, so after examining the polynomial to see that this is true to many decimal places, we round off. We then calculate the valuation of the norm of $\delta_i$ at the unique degree 1 prime in $K_8$ (using GP/PARI). For our choice of $\tau$, we find that $\delta_1 \notin \mathfrak{p}$, while $\delta_2 \in \mathfrak{p}$. Hence, $\mathrm{Tr}(\sigma(\tau)) = 6$. Then, using similar techniques, we determine that $\tau$ is a Frobenius element for the prime 7, but not for the primes 2, 13, or 23. Hence, for example, we predict that

$$\mathrm{Tr}(\rho(\mathrm{Frob}_2)) = \mathrm{Tr}(\sigma(\mathrm{Frob}_2))\epsilon_{31}(\mathrm{Frob}_2) = 4 \cdot (1) = 4,$$

and

$$\mathrm{Tr}(\rho(\mathrm{Frob}_{13})) = \mathrm{Tr}(\sigma(\mathrm{Frob}_{13}))\epsilon_{31}(\mathrm{Frob}_{13}) = 4 \cdot (-1) = 7.$$

| $\ell$ | 2 | 3 | 5 | 7 | 11 | 13 | 17 | 19 | 23 | 29 | 31 | 37 | 41 | 43 | 47 |
|---|---|---|---|---|---|---|---|---|---|---|---|---|---|---|---|
| $O(\rho(Fr_\ell))$ | 7 | 3 | 4 | 7 | * | 7 | 4 | 4 | 7 | 3 | * | 2 | 3 | 3 | 3 |
| $a(\ell,1)$ | 4 | 0 | 1 | 6 | 0 | 7 | 10 | 1 | 7 | 0 | * | 1 | 0 | 0 | 0 |
| $a(\ell,2)$ | 3 | 0 | 9 | 10 | 0 | 3 | 2 | 7 | 6 | 0 | * | 8 | 0 | 0 | 0 |

TABLE 10. Orders of $\rho(\mathrm{Frob}_\ell)$ and eigenvalues in weight $F(4,2,1)$

Returning to our study of $\rho$, we have that

$$\rho|_{I_p} \sim \begin{pmatrix} \theta^{190} & & \\ & \theta'^{190} & \\ & & \theta''^{190} \end{pmatrix}.$$

Note that $m = 190 = 3 + 6 * 11 + 1 * 11^2$. Hence, a weight predicted by the conjecture for $\rho$ is $F(6-2, 3-1, 1)' = F(4,2,1) = F(3,1,0) \otimes \det^1$. We may also take $m = 11 \cdot 190$ or $m = 11^2 \cdot 190$, which yield predicted weights of $F(6,6,0) \otimes \det^5$ and $F(8,3,0) \otimes \det^2$. Computing the cohomology in weight $F(3,1,0) \otimes \det^1$, we find a one-dimensional eigenspace with the eigenvalues indicated in Table 10. These eigenvalues are exactly what Conjecture 3.1 predicts, in order for $\rho$ to be attached to this eigenclass. The same system of eigenvalues (up to $\ell = 47$) also occurs in the other two weights predicted for $\rho$.

Similarly, the predicted weights for $\rho'$ are $F(6,0,0) \otimes \det^7$, $F(3,2,0) \otimes \det^6$ and $F(8,5,0) \otimes \det^8$. Each of these weights yields an eigenclass with the correct eigenvalues to have $\rho'$ attached (at least for $\ell \leq 47$).

## 8. COMPUTATIONAL TECHNIQUES

We now give an overview of our methods for computing the various Hecke eigenclasses on which we have reported in this paper. We begin by noting that we do not, in fact, do any direct calculations of cohomology. Instead we compute with homology, exploiting the natural duality, as in [7, Section 3].

Let $p$ and $N$ be positive integers with $N$ prime and let $V$ be a representation of the semigroup generated by $S_{pN}$ and $\Gamma_0(N)$. Then we wish to calculate $H_3(\Gamma_0(N), V)$, along with the action of various Hecke operators. The groups $H_3$ are easier to calculate than $H_1$ or $H_2$, since the virtual homological dimension of $\mathrm{SL}_3(\mathbb{Z}/p\mathbb{Z})$ is 3 (see [1]). In addition, one can show that for many classes of three-dimensional Galois representations, if the representation is attached to any homology class, then it is attached to a class in $H_3$ (c.f. [4]).

By Shapiro's lemma

$$H_3(\Gamma_0(N), V) \cong H_3(\mathrm{SL}_3(\mathbb{Z}), \mathrm{Ind}_H^G M),$$

and by [5, Lemma 1.1.4] this isomorphism is compatible with the action of the Hecke operators away from $pN$. This reduces our problem to computing the homology of the full group $\mathrm{SL}_3(\mathbb{Z})$ as long as we are willing to consider sufficiently general weights.

The broad outline of our calculations follows that of [1]. In particular we first use a slight modification of their Theorem 1 to identify $H_3(\mathrm{SL}_3(\mathbb{Z}), V)$ with the subspace of all $v \in V$ such that

1. $v \cdot d = v$ for all diagonal (but not necessarily scalar) matrices $d \in \mathrm{SL}_3(\mathbb{Z})$;
2. $v \cdot z = -v$ for all monomial matrices of order 2 in $\mathrm{SL}_3(\mathbb{Z})$;
3. $v + v \cdot h + v \cdot (h^2) = 0$,

where

$$h = \begin{pmatrix} 0 & -1 & 0 \\ 1 & -1 & 0 \\ 0 & 0 & 1 \end{pmatrix}.$$

We will refer to conditions 1 and 2 as the "semi-invariant condition" and to condition 3 as the "$h$-condition". Given a sufficiently concrete realization of $V$, computing the subspace satisfying these conditions is simply an exercise in linear algebra. In 8.2.2 we discuss some optimizations we have employed in carrying out the calculation. Once we have this subspace in hand, we then use Lemma 3 of [1] to compute the actions of various Hecke operators with respect to a basis of this space.

The main difference between our calculations and those in [1] is our use of more general coefficient modules. We will describe below our constructions of the modules $\mathrm{Ind}_{\Gamma_0(N)}^{\mathrm{SL}_3(\mathbb{Z})} F(a,b,c)$ for a $p$-reduced triple $(a,b,c)$. Another significant difference is a sharp increase in efficiency, and hence in the complexity of the calculations we can tackle. This increase is due partly to better algorithms (described below) and partly to having the entire calculation done using C++ code, rather than relying on Mathematica.

8.1. **Models for weights.** We have performed our calculations with a variety of weight modules. Our basic strategy has been to build more complicated weights up from simpler ones. In this subsection we will describe the $\mathrm{GL}_3(\mathbb{F}_p)$-modules with which we have worked, giving in particular a model for $F(a,b,c)$ for a general $p$-reduced triple $(a,b,c)$. Details of the implementation of these representations and of the process of inducing from $\Gamma_0(N)$ are left to subsection 8.2.

To begin we view $\bar{\mathbb{F}}_p^3$ as the standard 3-dimensional (right) $\bar{\mathbb{F}}_p[\mathrm{GL}_3(\mathbb{F}_p)]$-module on which $S_{pN}$ acts via reduction modulo $p$. Then $\mathrm{Sym}^g(\bar{\mathbb{F}}_p^3)$ is the space of homogeneous polynomials over $\bar{\mathbb{F}}_p$ of total degree $g$ in three variables $x, y, z$. An element $m$ of $GL_3(\mathbb{F}_p)$ acts on $f \in \mathrm{Sym}^g(\bar{\mathbb{F}}_p^3)$ by

$$f(\mathbf{x}) \cdot m = f(m\mathbf{x})$$

where $\mathbf{x}$ is the column vector ${}^t(x, y, z)$. Note that for $a \le p-1$ the representation $\mathrm{Sym}^a(\bar{\mathbb{F}}_p^3)$ is irreducible and is in fact isomorphic to $F(a,0,0) = W(a,0,0)$.

Next we look at the module $F(a,b,0)$ for $p$-restricted $(a,b,0)$.

**Theorem 8.1.** *Let $(a,b,0)$ be a $p$-restricted triple. Then the $\mathrm{GL}_3(\mathbb{F}_p)$-submodule of $\mathrm{Sym}^a(\bar{\mathbb{F}}_p^3) \otimes \mathrm{Sym}^b(\bar{\mathbb{F}}_p^3)$ generated by*

$$v = \sum_{i=0}^{b} (-1)^i \begin{pmatrix} b \\ i \end{pmatrix} (x^{a-i}y^i \otimes x^i y^{b-i})$$

*is isomorphic to $F(a,b,0)$.*

*Proof.* Recall that for any non-increasing triple $(\alpha,\beta,\gamma)$ of integers, both $W(\alpha,\beta,\gamma)$ and $F(\alpha,\beta,\gamma)$ are modules over $\mathrm{GL}_3(\bar{\mathbb{F}}_p)$ and not just over $\mathrm{GL}_3(\mathbb{F}_p)$. We will prove that the $\mathrm{GL}_3(\bar{\mathbb{F}}_p)$-module generated by $v$ is isomorphic to $F(a,b,0)$. Since $(a,b,0)$ is assumed to be $p$-restricted, $F(a,b,0)$ remains irreducible when viewed as a representation of $\mathrm{GL}_3(\mathbb{F}_p)$. We may then conclude that the $\mathrm{GL}_3(\mathbb{F}_p)$-module generated by $v$ is isomorphic to $F(a,b,0)$.

Since we are now looking at representations of $\mathrm{GL}_3$ of an algebraically closed field, we may employ the theory of highest weights in representations of algebraic groups [15, c.f. Section 31]. In particular, if we work with respect to the standard diagonal torus and the upper triangular Borel, we note that the non-increasing triples $(n_1,n_2,n_3)$ correspond to the dominant weights

$$\begin{pmatrix} t_1 & 0 & 0 \\ 0 & t_2 & 0 \\ 0 & 0 & t_3 \end{pmatrix} \mapsto t_1^{n_1} t_2^{n_2} t_3^{n_3}.$$

Then $F(n_1,n_2,n_3)$ is the unique irreducible representation of $\mathrm{GL}_3(\bar{\mathbb{F}}_p)$ with highest weight $(n_1,n_2,n_3)$.

Now, Young's rule [16, pg. 129] gives us that $W(a,0,0)\otimes W(b,0,0)$ has a filtration

$$W_0 \supset W_1 \supset \ldots \supset W_r = 0,$$

with the quotients $W_i/W_{i+1}$ isomorphic to the modules

$$W(a+b,0,0), \ldots, W(a+b-i,i,0), \ldots, W(a,b,0).$$

in the given order (so that $W(a+b,0,0)$ is a quotient and $W(a,b,0)$ is a submodule). Since $(a,b,0)$ is $p$-restricted, each $W(a+b-i,i,0)$ is irreducible if $a+b-i \le p-2$ or $a+b-2i = p-1$, and otherwise has $F(a+b-i,i,0)$ and $F(p-2,i,a+b-i-p+2)$ as composition factors [10, Prop. 2.11]. We see then that $F(a,b,0)$ appears only once as a composition factor of $W(a,0,0)\otimes W(b,0,0)$ and that it appears as a submodule and not just a subquotient.

It follows that $W(a,0,0)\otimes W(b,0,0)$ has a unique highest weight vector $w$ of weight $(a,b,0)$ and that the $\mathrm{GL}_3(\bar{\mathbb{F}}_p)$-module generated by this vector is isomorphic to $F(a,b,0)$. The lemma below shows that $v$ is such a vector, and hence the $\mathrm{GL}_3(\bar{\mathbb{F}}_p)$-module generated by $v$ is isomorphic to $F(a,b,0)$. □

**Lemma 8.2.** *The vector*

$$v = \sum_{i=0}^{b} (-1)^i \begin{pmatrix} b \\ i \end{pmatrix} (x^{a-i}y^i \otimes x^i y^{b-i})$$

*in* $\mathrm{Sym}^a(\bar{\mathbb{F}}_p^3)\otimes\mathrm{Sym}^b(\bar{\mathbb{F}}_p^3)$ *is a highest weight vector of weight* $(a,b,0)$*. Here "highest" refers to the usual lexicographic ordering of the weights.*

*Proof.* It is clear that $v$ is a weight vector of weight $(a,b,0)$. We need only show that the images of $v$ under the weight raising operators

$$g_1=\begin{pmatrix}1&0&0\\1&1&0\\0&0&1\end{pmatrix},\quad g_2=\begin{pmatrix}1&0&0\\0&1&0\\0&1&1\end{pmatrix},\quad g_3=\begin{pmatrix}1&0&0\\0&1&0\\1&0&1\end{pmatrix}$$

are all equal to $v$ plus something in the span $S$ of vectors of weight strictly less than $(a,b,0)$. Clearly $v\cdot g_2$ and $v\cdot g_3$ are both equal to $v$ modulo $S$. For $v\cdot g_1$, we calculate

$$\begin{aligned}
v\cdot g_1 &= \sum_{i=0}^{b}(-1)^i\binom{b}{i}x^{a-i}(x+y)^i\otimes x^i(x+y)^{(b-i)}\\
&= \sum_{i=0}^{b}(-1)^i\binom{b}{i}\sum_{k=0}^{i}\sum_{j=0}^{b-i}\binom{i}{k}\binom{b-i}{j}x^{a-i+k}y^{i-k}\otimes x^{i+j}y^{b-i-j}\\
&= \sum_{u=a-b}^{a}\sum_{v=a-u}^{b}\left(\sum_{i=a-u}^{v}(-1)^i\binom{b}{i}\binom{i}{u-a+i}\binom{b-i}{v-i}\right)x^uy^{a-u}\otimes x^vy^{b-v}
\end{aligned}$$

Setting $\alpha=i-(a-u)$, expanding the binomial coefficients and canceling equal terms, the inner sum becomes

$$\begin{aligned}
&\pm\sum_{\alpha=0}^{u+v-a}(-1)^\alpha\frac{b!}{(b-v)!(u+v-a-\alpha)!\alpha!(a-u)!}\\
=\ &\pm\frac{b!}{(b-v)!(a-u)!}\sum_{\alpha=0}^{u+v-a}(-1)^\alpha\frac{1}{\alpha!(u+v-a-\alpha)!}\\
=\ &\pm\frac{b!}{(b-v)!(a-u)!(u+v-a)!}\sum_{\alpha=0}^{u+v-a}(-1)^\alpha\binom{u+v-a}{\alpha}
\end{aligned}$$

which is zero if $u+v>a$. Thus the only terms $x^uy^{a-u}\otimes x^vy^{b-v}$ that appear in $v\cdot g_1$ with nonzero coefficient have $u+v=a$. It is now easy to see that $v\cdot g_1$ is in fact exactly equal to $v$. □

For arbitrary $p$-restricted $(a,b,c)$ we note that $F(a,b,c)\cong F(a-c,b-c,0)\otimes\det^c$. In practice we did all of our calculations with $F(a-c,b-c,0)$ and simply scaled by $\det^c$ at the end.

We have also made use of representations of the form $\mathrm{Sym}^a(\bar{\mathbb{F}}_p^3)\otimes\mathrm{Sym}^b(\bar{\mathbb{F}}_p^3)$, $\mathrm{Sym}^a(\bar{\mathbb{F}}_p^3)\otimes\mathrm{Sym}^b(\bar{\mathbb{F}}_p^3)^*$ and subquotients of $\mathrm{Sym}^a(\bar{\mathbb{F}}_p^3)$ for $a$ larger than $p-1$. By keeping track of the irreducible constituents of these representations we were sometimes able to show that certain systems of Hecke eigenvalues come from a specific irreducible module. See [12] for more details.

8.2. **Implementation.** The implementation of our algorithms has two very distinct parts. On the one hand we need to do calculations involving various $\mathrm{GL}_3(\mathbb{Z}/pN\mathbb{Z})$-modules $V$. This includes the basic vector space operations, as well as multiplying an element in $V$ by an element of $\mathrm{GL}_3(\mathbb{Z}/pN\mathbb{Z})$. Further we need to identify a basis of $V$ and be able to decompose elements of $V$ with respect to that basis. For efficiency reasons it is also important to be able to determine the coefficient of a given basis element in some product $v \cdot g$ without computing all of $v \cdot g$.

On the other hand, we need to carry out various higher level computations, such as finding the solutions to the $h$-condition above in order to compute homology with coefficients in $V$. These calculations can be described in terms of the basic operations of the previous paragraph without any specific knowledge about the module $V$. We have made use of object oriented programming techniques to keep these two computational issues strictly separated. This allows us to switch from computing with one module to another without having to rewrite any of the code describing the higher level algorithms.

8.2.1. *Coefficient modules.* We will now look at a few of the implementation details behind some of our coefficient modules. As we stated above, the basic building block for all of our representations is $\mathrm{Sym}^g(\bar{\mathbb{F}}_p^3)$, the space of homogeneous polynomials of degree g in three variables. The monomials form a natural basis of this space and it is a simple matter to compute the coefficient of any given monomial in a product $v \cdot g$. We have optimized this code to work especially well when many of $g$'s entries are zero. This is the case for the element $h$ above, as well as for many of the matrices arising in our Hecke operator calculations. The representations $\mathrm{Sym}^a(\bar{\mathbb{F}}_p^3) \otimes \mathrm{Sym}^b(\bar{\mathbb{F}}_p^3)$ again have natural bases coming from the monomial bases of $\mathrm{Sym}^a(\bar{\mathbb{F}}_p^3)$ and $\mathrm{Sym}^b(\bar{\mathbb{F}}_p^3)$ and all operations on the tensor product can be carried out in terms of those on each factor. We denote by $B_{ab} = \{w_i\}$ this basis of $\mathrm{Sym}^a(\bar{\mathbb{F}}_p^3) \otimes \mathrm{Sym}^b(\bar{\mathbb{F}}_p^3)$ and let $\langle \cdot, \cdot \rangle$ be the bilinear form with $\langle w_i, w_j \rangle = \delta_{ij}$.

The subspace $F(a, b, 0)$ of $\mathrm{Sym}^a \otimes \mathrm{Sym}^b$ does not come equipped with a canonical basis. For ease of computation we choose a basis in which each basis vector has a distinguished leading term. In other words, we choose a basis $\{v_i\}$ such that for each $i$ there is an element $w_i \in B_{ab}$ with $\langle w_i, v_i \rangle = 1$ and $\langle w_i, v_j \rangle = 0$ for $j \neq i$. We then let $\langle \cdot, \cdot \rangle_F$ be the bilinear form with $\langle v_i, v_j \rangle_F = \delta_{ij}$. Then for $v \in F(a, b, 0)$ we have

$$\langle v_j, v \rangle_F = \langle w_j, v \rangle$$

and so we can compute coordinates with respect to this basis of $F(a, b, 0)$ in terms of those with respect to the basis $B_{ab}$.

The final step in obtaining our general weights is to induce a representation $W$ from $\Gamma_0(N)$ to $\mathrm{SL}_3(\mathbb{Z})$. The $W$ we use are of the form $F(a, b, 0) \otimes \epsilon$ for some $\epsilon$ a character of $(\mathbb{Z}/N\mathbb{Z})^*$. We view $V = \mathrm{Ind}_{\Gamma_0(N)}^{\mathrm{SL}_3(\mathbb{Z})} W$ as the space of functions

$$V = \{f \colon \mathrm{SL}_3(\mathbb{Z}) \to W : f(xg) = f(x) \cdot g \text{ for } g \in \Gamma_0(N)\}$$

with $\mathrm{SL}_3(\mathbb{Z})$ acting by left translation.

We let $\{r_i\}$ be a set of representatives for $\mathrm{SL}_3(\mathbb{Z})/\Gamma_0(N)$, and $\{w_a\}$ be a basis for $W$. We again choose a bilinear form $\langle , \rangle$ on $W$ with $\langle w_a, w_b \rangle = \delta_{ab}$. Then we let $\phi_{r_i,w_a} \colon \mathrm{SL}_3(\mathbb{Z}) \to W$ by

$$\phi_{r_i,w_a}(x) = \begin{cases} w_a \cdot r_i^{-1} x & \text{if } x \in r_i \Gamma_0(N) \\ 0 & \text{otherwise .} \end{cases}$$

It is clear that the functions $\phi_{r_i,w_a}$ comprise a basis of $V$.

In order to express the action of $\mathrm{SL}_3(\mathbb{Z})$ on $V$ with respect to this basis, we need to introduce a bit of notation. For $x \in \mathrm{SL}_3(\mathbb{Z})$ let $\{x\}$ be the unique representative $r_i$ in $x\Gamma_0(N)$. Then

$$\begin{aligned}(\phi_{r_i,w_a}g)(x) &= \phi_{r_i,w_a}(gx)\\ &= \begin{cases} w_a \cdot r_i^{-1}gx & \text{if } gx \in r_i\Gamma_0(N)\\ 0 & \text{otherwise}\end{cases}\\ &= \begin{cases} w_a \cdot r_i^{-1}g\{g^{-1}r_i\}\{g^{-1}r_i\}^{-1}x & \text{if } x \in g^{-1}r_i\Gamma_0(N)\\ 0 & \text{otherwise}\end{cases}\\ &= \sum_b \left\langle w_a \cdot r_i^{-1}g\{g^{-1}r_i\}, w_b\right\rangle \phi_{\{g^{-1}r_i\},w_b}(x).\end{aligned}$$

Note that in order to compute the actions of Hecke operators on $H_3(\mathrm{SL}_3(\mathbb{Z}), V)$ we also need to know how elements of

$$S = \{m \in M_3(\mathbb{Z}) : \det(m) \text{ is positive and prime to } pN\}$$

act on $V$. Let $\Sigma$ be the semigroup generated by $\Gamma_0(N)$ and $S_{pN}$. Then $S = \mathrm{SL}_3(\mathbb{Z})\Sigma$ and $\Gamma_0(N) = \mathrm{SL}_3(\mathbb{Z}) \cap \Sigma$ (this is part of what it means for the Hecke pair $(\Gamma_0(N), \Sigma)$ to be compatible to $(\mathrm{SL}_3(\mathbb{Z}), S)$). Thus if $m \in S$, we have $m = ns$ for some $n \in \mathrm{SL}_3(\mathbb{Z})$ and $s \in \Sigma$. Moreover, $n$ is determined modulo $\Gamma_0(N)$ and so the coset representative $\{n\}$ depends only on $m$. If we extend our notation to write $\{m\} = \{n\}$, the formula above for the action of $g$ on $\phi_{r_i,w_a}$ makes sense for any $g \in S$. This action of $S$ on $V$ described by the formula induces the correct action of $\mathcal{H}(pN)$ on $H_3(\mathrm{SL}_3(\mathbb{Z}), V)$ (i.e. the one compatible with the action on $H_3(\Gamma_0(N), W)$.

The $r_i$ may be chosen so that each is congruent to the identity modulo $p$, which greatly speeds up some of the calculations. Note that $\mathrm{SL}_3(\mathbb{Z})/\Gamma_0(N) \cong \mathbb{P}^2(\mathbb{Z}/N)$ and so is easy to work with. Also note that our formula shows at once how to compute the coefficients of a basis element $\phi_{r_j,w_b}$ in $v \cdot g$ for $v \in V$ and $g \in S$.

8.2.2. *Finding homology.* Now we move on to the general algorithms we have used to compute the homology of $\mathrm{SL}_3(\mathbb{Z})$ with coefficients in some representation $V$. While this is a simple exercise in linear algebra, we have found it useful to tailor certain optimizations to our situation to allow us to work with larger examples. A typical instance of finding the solutions to the $h$-condition, for instance, involves finding the kernel of a $700000 \times 30000$ matrix. These optimizations have been largely heuristic. We make no claim of having optimal algorithms.

Let $V$ be a $\Sigma$-module of dimension $d$, with basis $\{v_i\}$ and let $\langle\cdot,\cdot\rangle$ be the bilinear form with $\langle v_i, v_j\rangle = \delta_{ij}$. Let $M$ be the 24-element group of monomial matrices in $\mathrm{SL}_3(\mathbb{Z})$. Then for $p > 3$ the space of semi-invariants in $V$ is the image of the operator

$$P = \sum_{g\in M} \epsilon(g)g$$

where $\epsilon(g)$ is the sign of the permutation on three letters induced by $g$. Our computations do not include examples for which $p = 2$, and for $p = 3$ only a minor adjustment is needed. Computing the action of $P$ on each $v_i$ is not computationally intensive since we have specially optimized all of our coefficient modules with regard to the operation of monomial matrices. We then use column reduction to find a basis for $V \cdot P$. We note that the dimension $d_{\text{semi}}$ of $V \cdot P$ is approximately $d/24$.

The more serious stage of the calculation is finding the solutions of the $h$-condition on $V \cdot P$. We describe our algorithm for finding the solutions of the $h$-condition on any $c$-dimensional subspace $W$ of $V \cdot P$ with basis $\{b_i\}$. We are looking for the nullspace of the $d \times c$ matrix $M = (m_{ij})$ where $m_{ij} = \langle v_i, b_j \cdot (1 + h + h^2) \rangle$ is the coordinate of $v_i$ in $b_j \cdot (1 + h + h^2)$. Simply computing this matrix and performing Gaussian elimination would theoretically allow us to find the nullspace but is hopelessly inefficient in both space and time. Although the matrix $M$ is quite sparse, it becomes much denser as the elimination proceeds. Since we work with very large $d$ ($d$ on the order of $7 \times 10^5$ is not uncommon), we would rapidly run out of memory. We will touch on four optimizations we have made to speed up the calculation and to reduce the memory requirements.

First, we note that the rows of $M$ are highly redundant as there are at most about $1/24^{th}$ as many columns as rows. We exploit this by computing the rows of $M$ one at a time and only storing those that yield new information about the kernel. Recall that we have set up our coefficient modules so that we can individually compute the entries $\langle v_i, b_j \cdot (1 + h + h^2) \rangle$ in the $i^{\text{th}}$ row of $M$ *without* having to compute all of $b_j \cdot (1 + h + h^2)$. We discuss below another optimization that makes this separate computation especially efficient. As we find a new row $R$ we continue our elimination process by subtracting from $R$ the appropriate multiples of all the previously stored rows. If we are left with a non-zero row, we append it to our stored matrix, which remains in row-echelon form. If we are left with the 0 row, then $R$ did not add any constraints on the kernel of $M$ and we may discard it and move on to the next row. This guarantees that we never waste space by storing redundant rows, and caps the maximum number of rows we will ever store at $c \leq d_{\text{semi}} \approx \frac{d}{24}$. We will denote by $E$ the matrix that we are building up row by row in this process.

Our second optimization is motivated by the fact that most of the information about the kernel of $M$ can be obtained from $M$'s early rows. At each stage in our calculation, we clearly have $\ker M \subset \ker E$. Since $E$ is in row-echelon form we can immediately read off the dimension of $\ker E$. In practice we find that the dimension of $\ker E$ drops below 1 or 2 percent of $d_{semi}$ after we run through as few as $1/5^{\text{th}}$ or so of the rows of $M$. Once this happens, we pause our calculation and compute (a basis for) the kernel of $E$, which is relatively easy to do since $E$ is already in row-echelon form. We have now reduced our problem to finding the kernel of $1 + h + h^2$ not on $W$ but on the much smaller space $\ker E$. We then start the algorithm over replacing $W$ by $\ker E$. Our new choice of $W$ guarantees that the initial rows of the new matrix $M$ will all be 0, and so we can resume our calculation with the row at which we had paused. It is crucial here that we have not computed $M$ all at once and thus do not have to make any time consuming adjustments to account for our new basis. Indeed, it is now much easier to compute the new $M$, as it has far fewer rows. In practice the calculation very rapidly proceeds through the remaining rows of $M$ and then computes the kernel of the new $E$, which is equal to the kernel of $M$. Our choice of a cutoff on the dimension of $\ker E$ is entirely heuristic, and we adjust it based on the size of $V$.

Both of the optimizations above rely on the efficiency of the calculation of the coefficients $\langle v_i, v \cdot (1 + h + h^2) \rangle$ of each $v_i$ in $v \cdot (1+h+h^2)$ for $v \in V$. Although our modules allow for the calculation of $\langle v_i, v \cdot g \rangle$ for any $v$ and $g$ without computing all of $v \cdot g$, there is still a great deal of work duplicated if we separately perform this calculation for all of the $v_i$. For our calculations of the Hecke operators (see below) this is not necessary, but as described above we must do this in the cases $g = h$ and $g = h^2$. We have optimized for this by storing some of the common pieces of these calculations. For example, when $V = \mathrm{Ind}_{\Gamma_0(N)}^{\mathrm{SL}_3(\mathbb{Z})} W$, we begin by computing and storing the entire matrices describing the actions of $h$ and $h^2$ on $W$, and also the permutations induced by $h$ and $h^2$ on $\mathbb{P}^2(\mathbb{Z}/N)$. Since the dimension of $W$ is small compared to the dimension of $V$ (even when $N = 2$, the dimension of $V$ is 7 times that of $W$), this calculation is not terribly costly in space or time. These stored tables can then be used to compute the action of $h$ and $h^2$ on elements of $V$ very quickly. We have implemented similar strategies when $V$ is not induced, but is the tensor product of two smaller representations.

Finally, we have increased our available memory by making use of disk space and swapping pieces of our matrix in and out of memory. This requires minor modifications to the reduction algorithm described above in order to reduce the number of disk swaps. In particular, we carry out our row reduction on several (1000) new rows at once. In the end, this does not have a dramatic effect on run time, but slashes the amount of RAM required.

8.2.3. *Computing the Hecke Action.* Our computation of the action of the Hecke operators closely mirrors that in [1] and we refer the reader to Sections 3 and 8 there for a discussion of modular symbols and a description of the action of Hecke operators on homology. We will just summarize by noting that for $v \in V$ satisfying the semi-invariant condition and the $h$-condition, we have

$$T(l,k)v = \sum_{i,j} v \cdot M_{ij} B_i,$$

where

$$\Gamma_0(N) D(l,k) \Gamma_0(N) = \coprod_i \Gamma_0(N) B_i$$

and the modular symbol $[B_i]$ is homologous to $\sum_j [M_{ij}]$. We have not recomputed the matrices $M_{ij}$ but have used the files generated in the course of the calculations in [1].

If $\{f_l\}$ is a basis for the semi-invariants in $V$ satisfying the $h$-condition, then we know a priori that $\sum_{ij} f_k \cdot M_{ij} B_i$ will be a linear combination $\sum_l a_{kl} f_l$ of the $f_l$. We wish to obtain the numbers $a_{kl}$. To do this efficiently we use the same trick we employed in our choice of basis for $F(a,b,0)$ and adjust our basis $\{f_l\}$ so that for each $l$ there is a basis vector $v_l$ of $V$ so that $\langle v_l, f_m \rangle = \delta_{lm}$. Then $a_{kl}$ is the sum over $i$ and $j$ of $\langle v_l, f_k \cdot M_{ij} B_i \rangle$. As we have discussed, we are able to compute these coefficients directly. This is vastly superior to computing all of $f_k \cdot M_{ij} B_i$ since the dimension of the homology space is only a tiny fraction of the dimension of $V$. This technique was used in [1], although it could not be implemented as efficiently there due to the reliance on Mathematica's multivariate polynomial routines.

A final optimization uses the fact that the Hecke operators we are dealing with all commute and so preserve each other's eigenspaces. The ultimate goal of our calculation is to identify simultaneous eigenvectors $v$ of the $T(l,k)$ attached to given Galois representations, that is with $T(l,k)v = \alpha(l,k)$ for some prescribed $\alpha(l,k)$. If we compute the entire matrix for the action of $T(2,1)$ (which is very easy since $T(2,1)$ involves only 13 $M_{ij}$ terms, whereas $T(47,1)$ involves 55923 such terms) and find a single eigenvector $v$ with eigenvalue $\alpha(2,1)$ we need only compute the image of the other $T(l,k)$ on $v$ and not on the whole homology space. Moreover, we know that $v$ will be an eigenvector of each $T(l,k)$ and so we only need to compute a *single* coefficient $\langle v_\ell, T(l,k)\rangle$ in order to determine the eigenvalue. This gives an extraordinary reduction in the time required to make the calculation. For example, we find that dimension of the homology space at level $\Gamma_0(11)$, weight $F(22,11)(\epsilon_{11})$ and $p = 19$ is 31. We are interested in a single eigenvector in this space. In order to compute the entire matrix of a Hecke operator we would need to find $31^2 = 961$ coefficients of basis vectors. Instead, we reduce this to a single coefficient, giving nearly a thousand fold increase in performance. We point out that this technique was not needed in [1] as the homology spaces dealt with there were much smaller.

8.2.4. *Reliability.* Whenever relying on a large amount of computer calculation, one hopes for a number of consistency checks on the data. Our first check is that two entirely independent programs were written to carry out the calculations on several different computers by two different authors and both programs yielded identical data where compared. The programming was done in C and C++ and compiled with gcc running on a Sparc Ultra 5, a Pentium II under Linux, and a Pentium III under Linux. We also compared our data to some of the data obtained in [1] and [4] and found everything to be consistent.

Other checks include the fact that, whenever tested, the operators $T(l,k)$ on a given homology space all commuted and that (again when tested) the Hecke operators all did preserve the space of semi-invariants in $V$ satisfying the $h$-condition. Perhaps more compelling is the fact that our data meshes exactly with the Galois representations we have studied. While the correspondence is only conjectural, the agreement we observed very strongly suggests the validity of our calculations.

## 9. Acknowledgments

The authors thank Warren Sinnott and Richard Taylor for many helpful conversations in the course of this work. We thank the maintainers of the Bordeaux database of number fields (available at `ftp://megrez.math.u-bordeaux.fr`), from which we have obtained defining polynomials for many of the fields described in this paper. We also thank John Jones and David Roberts, whose tables of number fields ramified only at small primes (available on the web at `http://math.la.asu.edu/~jj/numberfields/` also provided many of our defining polynomials.

(Avner Ash) Department of Mathematics, Boston College, Chestnut Hill, MA 02467
*E-mail address*: `avner.ash@bc.edu`

(Darrin Doud) Department of Mathematics, Harvard University, Cambridge, MA 02138
*E-mail address*: `doud@math.harvard.edu`

(David Pollack) Department of Mathematics, The Ohio State University, Columbus, OH 43210
*E-mail address*: `pollack@math.ohio-state.edu`